\documentclass{amsart}
\usepackage{amsmath,amsthm,amssymb,amsfonts,amscd}

\newcommand \red   {\ensuremath{\mathrm{red}}}

\newcommand \ini {\ensuremath{\mathrm{in}}}

\newcommand \Tor {\ensuremath{\mathrm{Tor}}}
\newcommand \Gin {\ensuremath{\mathrm{Gin}}}
\newcommand \sat {\ensuremath{\mathrm{sat}}}
\newcommand \depth {\ensuremath{\mathrm{depth }}}
\newcommand \Spor {\ensuremath{\mathrm{Spor}}}
\newcommand \codim {\ensuremath{\mathrm{codim }}}
\newcommand \supp {\ensuremath{\mathrm{supp}}}

\newcommand \bbp {\mathbb P}
\newcommand \reg {\mathrm {reg}}

\newcommand \kxn {K[x_1,\ldots,x_n]}

\theoremstyle{plain} 
\newtheorem{Thm}{Theorem}[section]
\newtheorem{Prop}[Thm]{Proposition}
\newtheorem{Coro}[Thm]{Corollary}
\newtheorem{Lem}[Thm]{Lemma}

\theoremstyle{definition}
\newtheorem{Def}[Thm]{Definition}
\newtheorem{Ex}[Thm]{Example}
\newtheorem{Remk}[Thm]{Remark}

\numberwithin{equation}{section}

\begin{document}

\title{Some Geometric Results arising from the Borel Fixed Property  }
\author[J.\ Ahn and J.\ Migliore]{Jeaman Ahn$^{*}$, Juan C.\ Migliore$^{**}$}

\thanks{$^*$ The first author expresses his gratitude to the University of
Notre Dame Mathematics Department and to Juan Migliore for their
hospitality.  He
also thanks BK21 and Younghyun Cho, his advisor; without their
support, this paper
would not have been written.
\\
$^{**}$ Part of the work for this  paper was done while the second author was
sponsored by the National Security Agency  under Grant Number
MDA904-03-1-0071.
}

\begin{abstract}
In this paper, we will give some geometric results using generic
initial ideals for the degree reverse lex order.  The first
application is to the
regularity of a Cohen-Macaulay algebra, and we improve a well-known bound.  The
main goal of the paper, though, is
to improve on results of Bigatti, Geramita and Migliore concerning geometric
consequences of maximal growth of the Hilbert function of the
Artinian reduction of
a set of points.  When the points have the Uniform Position Property, the
consequences they gave are even more striking.  Here we weaken the growth
condition, assuming only that  the values of the  Hilbert function of
the Artinian
reduction are equal in two consecutive degrees, and that the first of
these degrees
is greater than the second reduction number of the points.  We
continue to get nice
geometric consequences even from this weaker assumption.  However, we have
surprising examples to show that imposing the Uniform Position Property on the
points does not give the striking consequences that one might expect.
This leads
to a better understanding of the Hilbert function, and the ideal
itself, of a set
of points with the Uniform Position Property, which is an important
open question.  In the last section we describe the role played by the Weak
Lefschetz Property (WLP) in this theory, and we show that the general hyperplane
section of a smooth curve may not have WLP.
\end{abstract}

\maketitle


\section{Introduction}
It has been shown in work of Conca, Fl\o ystad, Green and many
others that generic initial ideals contain a tremendous amount of
geometric information about the subschemes from which they are
obtained.  This paper continues in this vein, applying the theory
of generic initial ideals to the study of geometric consequences
that arise when the Hilbert function has certain growth properties
(but possibly short of maximal growth). Along the way, we also
obtain results about the regularity of a Cohen-Macaulay ring.  Our
main application is to the case of (reduced) points in projective
space, especially those having the Uniform Position Property (UPP)
(i.e.\ the  property that any two subsets of the same cardinality
have the same Hilbert function).  Here what is {\em not} true
turns out to be as striking as what {\em is}  true, and we give
some new insight into the behavior of points with UPP.

We begin in section \ref{prelim sect} with a review of the basic facts
about generic initial ideals.  In the process, we make some new
observations.  Our main objects of interest are two natural invariants of
a Borel fixed monomial ideal (see the definition of $D(A)$ and $M(A)$ on
page \pageref{def of D and M}).  We show how these invariants are related
to many other important invariants such as the depth, codimension, etc.
We also obtain results about the vanishing and non-vanishing of
cohomology of the ideal sheaf of a closed subscheme of projective space,
and a criterion about the Cohen-Macaulayness of the coordinate ring of a
homogeneous ideal.  In this section we also recall the definition of an
invariant which will be central to the rest of the paper, namely the {\em
$s$-reduction number} of $R/I$, which was introduced by Hoa and Trung
\cite{HT}.  This has several equivalent formulations (Lemma
\ref{Reduction Number}
gives some of them), but the most convenient for us is that $r_s(R/I)
= \max \{ k \
| \ H(R/(I+J),k)
\neq 0 \}$, where $J$ is the ideal generated by $s$ generally chosen
linear forms.

The computation of the second reduction number, for a
zero-dimensional scheme $Z$,
is closely connected to the Weak Lefschetz Property (WLP) for the Artinian
reduction of $R/I_Z$.  WLP may be defined as follows:  If $J = (L_1,L_2)$ is
an ideal of generally chosen linear forms, then WLP is equivalent to the
property that
\[
H(R/(I_Z+J),t) = \max \{ \Delta^2 (R/I_Z,t), 0 \}.
\]
In particular, if WLP holds for the Artinian reduction of $R/I_Z$ then
\[
r_2(R/I_Z) = \max \{ t | \Delta H(R/I_Z,t) > \Delta
H(R/I_Z,t-1)\}.
\]
   It is very
natural to ask whether a set of points with UPP automatically has
WLP, since this
property is the ``expected" one and UPP is somehow a ``general" property.  We
remark that we were somewhat surprised by Example \ref{UPP vs WLP},
which gives a
set of points with UPP that does not have WLP.  Furthermore, since a general
hyperplane or hypersurface section of a smooth curve has UPP, it is also very
natural to ask whether the same is true for WLP.  We were equally surprised to
see that the same example can be modified to answer both of these questions in
the negative; this is shown in Example \ref{hyper sect}.

In section \ref{reg sect} we apply the ideas of section \ref{prelim
sect} to give a bound on the regularity of an arithmetically
Cohen-Macaulay  projective subscheme.  We generalize the well known fact
(cf.\ \cite{E2}) that
\[
\reg(I_X)\leq \deg(X)-\codim(X)+1
\]
by proving that if $X$ is an arithmetically Cohen-Macaulay subscheme of
$\bbp^{n-1}$ with $\codim (X)=e$ then
\[
\reg (I_X)\leq\deg (X) -\binom{\alpha-1+e}{\alpha-1}+\alpha,
\]
where $\alpha=\alpha(I_X)$ is the initial degree of $I_X$ (see Theorem
\ref{Cohen1}).  This was also proved by Nagel (who in fact gave a more general
statement, but with a very different proof) in \cite{uwe}.

In the paper \cite{BGM}, the authors studied some geometric consequences
that arise from maximal growth of the first difference, $\Delta H$, of the
Hilbert function of some projective subscheme.  This has seen the most
interest and applications when the subscheme is a finite set of (reduced)
points, $Z$ (so $\Delta H$ is the $h$-vector, i.e.\ the Hilbert function
of the Artinian reduction), and here the most interesting case (from our
point of view)  comes when the maximal growth is given by $\Delta
H(R/Z,d) = \Delta H(R/I_Z,d+1) = s$ (say).  In this case, to be maximal
growth means that necessarily $d \geq s$.  When this equality on $\Delta
H$ occurs for $d \geq s$, it was shown in \cite{BGM} that the component
$(I_Z)_d$ defines a {\em reduced} curve, $V$, of degree $s$,  and that in
fact $\reg(I_V) \leq d$.  In this context it is not necessarily true that
$V$ is unmixed.  Let $C$ be the pure one-dimensional part of $V$.  When, in
addition, we assume that
$Z$ has UPP, it was also shown that $Z \subset C$, $(I_Z)_d =
(I_C)_d$, and that
$C$ is reduced and irreducible.  This has strong implications for the possible
Hilbert functions of sets of points with UPP, although it is still far from a
classification.  For example, it was not proved in \cite{BGM}, but we remark in
Proposition \ref{strictly decreasing} that these conditions imply
that once $\Delta
H(R/I_Z,t) < \Delta H(R/I_Z,t-1)$ for some $t > d+1$, then the
Hilbert function is
strictly decreasing from then on.

Section \ref{first points sect} begins the main geometric task in this
paper, which is weaken the assumption $d \leq s$ of \cite{BGM}.
  Our results are given in terms of the second
reduction number, $r_2(R/I_Z)$.  As mentioned above, if $J$ is the ideal
generated by two general linear forms, then $r_2(R/I_Z)$ is the maximum
degree in which $R/(I_Z +J)$ is non-zero (we do not necessarily have
WLP here).  We
note in Remark \ref{remark1} that $d > r_2(R/I_Z)$ really is a weaker
hypothesis
than $d \geq s$.  We also note in Remark \ref{not curve} that in
general it is much
weaker, and that the hypothesis cannot be weakened further.

Our first main result is that if
\begin{equation} \label{condition}
\Delta H(R/Z,d) = \Delta H(R/I_Z,d+1) =
s \hbox{ for some
$d > r_2(R/I_Z)$}
\end{equation}
  then again (as in \cite{BGM}) $(I_Z)_d$ is the degree $d$ component of
the saturated ideal of some curve $V$ (not necessarily unmixed) of
degree $s$, and
$\reg (I_V) \leq d$.   Easy examples show that this cannot be
extended to $d \leq
r_2(R/I_Z)$.

Enboldened by this result, in section \ref{UPP sect} we sought to prove
results analogous to those of \cite{BGM} for the case where $Z$ has, in
addition, UPP.  We found, on the contrary, some very surprising examples,
and in doing so we have gained some new insight (albeit negative) into the
behavior of the Hilbert function of points with UPP.  It is not as
well-behaved as one might think.  Luckily, we do retain the fact that if
(\ref{condition}) holds for a set, $Z$, with UPP then $Z \subset C$
(Theorem \ref{UPP}), where $C$ is the unmixed part of $V$.  That is, $V = C$ is
unmixed.

It is known that if $Z$ is a set of points with UPP, a general element of
smallest degree of $I_Z$ is reduced and irreducible.  For points
in ${\mathbb P}^3$, for instance, this means that a general surface of
least degree containing $Z$ is reduced and irreducible.  In particular,
when this surface is unique, it is reduced and irreducible.  We expected
that if $Z$ is a finite set of points with UPP, then in the least
degree component
of $I_Z$ where the base locus is a curve, $C$, we would analogously
get that $C$ is
reduced and irreducible; it was shown in \cite{BGM} that this is true when $d
\geq s$.  On the contrary, however,  Example
\ref{CD gen} shows that if we assume only condition (\ref{condition})
and UPP, we do
get an unmixed curve $C$ containing $Z$ (as noted above), but $C$ can fail
to be either reduced or irreducible.  We also expected that even in this
situation it would continue to be true that if
$Z$ has UPP and satisfies (\ref{condition}), then once
$\Delta H(R/I_Z,t) < \Delta H(R/I_Z,t-1)$ for some $t > d+1$, then the first
difference of the Hilbert function would be strictly decreasing from then on. 
But on the contrary, Example \ref{CD gen} also shows that this is false.  (As
mentioned above, though, Proposition \ref{strictly decreasing} gives conditions
where the first difference of the Hilbert function {\em is} strictly
decreasing.)  One might
think that somehow Example
\ref{CD gen} was an ``accident," in the sense that in Example \ref{CD
gen}, while
the curve
$C$ (defined by the degree
$d$ component of $Z$) does fail to be either reduced or irreducible,
still there is
a reduced and irreducible component of $C$ that contains all of the
points.  Perhaps
at least this is always true for points with UPP.  But in fact, Example
\ref{not irred} removes even this hope. This is an example of a set of points
with UPP, satisfying condition (\ref{condition}), for which the curve
defined by the component of degree $d$ is reduced and consists of two
``identical" components, each of which contains exactly half of the
points.

Our results in section \ref{add'l comments} were motivated both by the desire to
obtain some results in the range $d \leq r_2(R/I_Z)$, and also by the desire to
understand the role of WLP in this theory beyond simply being a useful tool for
computing $r_2(R/I_Z)$.  Now it turns out that the {\em second} difference of
the Hilbert function and $r_3(R/I_Z)$ are important.  One has to be more careful
here, since reducing by a second linear form does not necessarily have the same
Hilbert function as the second difference of the original Hilbert function; this
does hold when we have WLP, and we use this fact.  For instance, we show in
Theorem \ref{WLP delta 2} that if
$Z$ is a zero-dimensional subscheme of ${\mathbb P}^{n-1}$, $n>3$, with WLP and
if 
\[
\Delta^2 H(R/I_Z,d) = \Delta^2(R/I_Z,d+1) = s
\]
for $r_2(R/I_Z) > d > r_3(R/I_Z)$, then $\langle (I_Z)_{\leq d} \rangle$ is a
saturated ideal defining a two-dimensional  scheme of degree $s$ in ${\mathbb
P}^{n-1}$, and it is $d$-regular.  We also show in Corollary \ref{UPP_WLP} that
if $Z$ has both UPP and WLP, and if $I_Z$ has generators in degrees $d_1 \leq
d_2 \leq \dots$, then
\[
\Delta^2 H(R/I_Z,d) > \Delta^2 H(R/I_Z,d+1)
\]
for $d_2 \leq d < r_2(R/I_Z)$.  As mentioned above, we also show (Example
\ref{UPP vs WLP}) that UPP does not necessarily imply WLP, and in Example
\ref{hyper sect} we modify this example to show   that the general hyperplane or
hypersurface section of a smooth curve does not necessarily have WLP (although
it does necessarily have WLP if the curve is in ${\mathbb P}^3$).

Finally, we would like to comment on this invariant $r_2(R/I_Z)$.  It is seen in
this paper to be a very natural invariant, and when WLP holds it is easy to
compute.  WLP has been shown in other papers to be a very commonly occuring
property, as has UPP, and in section \ref{add'l comments} we have shown that it
plays an integral role in our study of Hilbert functions.  However, for an
arbitrary set of points it is not always easy to determine whether WLP holds,
nor is it always easy to compute
$r_2(R/I_Z)$. This makes our results less easy to apply than, for instance, the
corresponding results of \cite{BGM} where it is assumed  that $d \geq s$. 
However, it is not hard to see that $r_2(R/I_Z)$ can be much smaller than $s$. 
Our methods thus allow us to give results in a (possibly large) range in which
the methods of
\cite{BGM} absolutely do not apply.  This results
in new information about the Hilbert function of a set of points with UPP.  We
also show that very interesting subtle differences arise between what happens in
this range and what happens in the range $d \geq s$, leading to interesting
examples of unusual behavior in the Hilbert function of a set of points with
UPP.  We believe that our methods will lead to further progress on these
questions.

We are grateful to Uwe Nagel for asking us the question that led to
Example~\ref{hyper sect}.


\section{Preliminaries} \label{prelim sect}
In this section, we survey some definitions, notation and some
preliminary facts for generic initial ideals. Let
$R=K[x_1,\ldots,x_n]$ be the polynomial ring over a field $K$ of
characteristic $0$. For any $g=(g_{ij}) \in \textup{GL}_n(R_1^\vee)$,
we define an action on $R$ which induces an $k$-algebra
isomorphism for any homogeneous form $f \in R$ by
\[
f(x_1,\ldots,x_n) \mapsto f(g(x_1),\ldots,g(x_n)),
\]
where $g(x_i)= \sum_{j=1}^{n} g_{ij}x_j$. A monomial ideal $I$ is
said to be {\em Borel fixed} if
\[g(I)=I\]
for every upper triangular matrix $g\in \textup{GL}_n(R_1^\vee)$.
A Borel fixed monomial ideal has a nice combinatorial property, namely that
of being {\em strongly stable}: If $x_im \in I$ for some
monomial $m\in I$, then $x_jm\in I$ for all $j<i$.

For any monomial term order $\tau$, the initial ideal of a
homogeneous ideal $I\subset R$ depends on the choice of variables
and basis made. By allowing a generic change of basis and
coordinates, we may eliminate this dependence.
\begin{Thm}{\textup{(\cite{Ga}, \cite{BS}, \cite{E1})}} For any monomial term
order $\tau$ and any homogeneous ideal $I$, there is a Zariski
open subset $U\subset \textup{GL}_n(R_1^\vee)$ such that
$\ini_\tau(g(I))$ is constant and Borel fixed for $g\in U$. We
will call $\ini_\tau(g(I))$ the generic initial ideal of $I$ and
denote it $\Gin_\tau(I)$.
\end{Thm}

Note that a homogeneous ideal $I$ is {\em saturated} if
\[(I:(x_1,\ldots,x_n))=I.\]
The {\em saturation} of $I$ is
\[I^{sat}=\bigcup_{k\geq 0} (I:(x_1,\ldots,x_n)^k).\]
A homogeneous ideal $I$ is {\em $m$-saturated} if
\[I^{sat}_d=I_d \quad \textup{for all}\,d\geq m.\]
The {\em saturation degree} of $I$, denoted $\sat(I)$, is the smallest
$m$ for which $I$ is $m$-saturated.

A homogeneous ideal $I$ is {\em $m$-regular} if, in the minimal free
resolution of $I$, for all $p\geq 0$, every $p$-th syzygy has
degree $\leq m+p$. The {\em regularity} of $I$, $\reg(I)$, is the smallest such
$m$.  The regularity has an alternate description, in terms of
cohomology.  David
Mumford defined the regularity of a coherent sheaf on projective
space (now known as
Castelnuovo-Mumford regularity) as follows: a coherent sheaf
$\mathcal F$ on $\mathbb P^{n-1}$ is said to be {\em $m$-regular} if
$H^q(\mathbb
P^{n-1},\mathcal F(m-q))=0$ for all $q>0$;  the regularity,
$\reg({\mathcal F})$, is
the smallest such $m$. If $I$ is a saturated ideal, $m$-regularity in the first
sense is equivalent to the geometric condition that the associated
sheaf $\mathcal
I$, on projective space $\mathbb P^{n-1}$, satisfies the condition of
Casteluovo-Mumford $m$-regularity. In general case, we may show the
following fact
with local cohomology \cite{E2}:
\[
\reg(I)=\max\{\sat(I), \reg(\mathcal I)\}.
\]

A great deal of fundamental information about $I$ can be read off  if
$I$ is a Borel
fixed monomial ideal of $R$. The following is a useful property of Borel fixed
ideals.
\begin{Thm}\textup{(\cite{BS}, \cite{G})}\label{BS}
For a Borel fixed monomial ideal $I\subset R=\kxn$,
\begin{enumerate}
\item[(a)] $I^{\sat}=\bigcup_{k=0}^{\infty}(I:x_n^k)$.
\item[(b)]
$\sat(I)=$ the maximal degree of generators involving $x_n$.
\item[(c)] $\reg(I)=$ the maximal degree of generators of $I$.
\end{enumerate}
\end{Thm}

We introduce a piece of notation. If $K=(k_1,\ldots,k_n)$, we
denote by $x^K$ the monomial
\[
x^K = x_1^{k_1}\cdots x_n^{k_n}
\]
in $R$, and by $|K|$ its degree $|K|=\sum_{j=1}^{n}k_j$. For monomial
$x^K$ with exponent $K=(k_1,\ldots,k_n)$ one defines
\[
\max(x^K)=\max\{j : k_j>0\}
\]
\[
\min(x^K)=\min\{j : k_j>0\}.
\]
For a set of monomials $A$ \label{def of D and M}
\[
D(A)=\max\{\min(x^K) : x^K \in A\}
\]
\[
M(A)=\max\{\max(x^K) : x^K \in A\}
\]
If $A$ is the set of the minimal generators of a monomial ideal
$I\subset R$, we set
\[
D(I)= D(A), \quad M(I)=M(A).
\]
Then we get the following lemma.

\begin{Lem}{\label{BR^D(I)}}
If $I$ is a Borel fixed monomial ideal of $R=k[x_1,\ldots,x_n]$
such that $n\geq 2$ and $\dim (R/I)>0$ then,
\begin{enumerate}
\item[(a)] $D(I)=\codim (R/I)=n-\dim(R/I)$
\item[(b)]
$M(I)=\textup{codepth} (R/I)=n-\depth (R/I)$.
\end{enumerate}
\end{Lem}

\begin{proof}
By definition of $D(I)$, there is a minimal generator $x^K \in A$
such that $D(I)=\min(x^K)$, where $A$ is the set of the minimal
monomial generators of $I$. Since $I$ is strongly stable,
$x_{D(I)}^{|K|}$ must be in $I$ and thus
$\sqrt{I}=(x_1,x_2,\cdots,x_{D(I)})$. Since the dimension of $R/I$
is equal to the dimension of $R/\sqrt{I}$, where $\sqrt{I}$ is the
radical ideal of $I$,
\[\dim(R/I)=n-D(I).\]

For the proof of (b), We begin with induction on the number of
variables $n$. In case $n=2$, it is clear. Suppose that $n>2$. If
$M(I)=n$ then there is a minimal generator of $I$ involving $x_n$,
by definition of $M(I)$. Hence $I$ is not a saturated ideal. Let
$H^0_m(R/I)$ be the $0$-th local cohomology of $R/I$ and let
$I^{\sat}$ be the saturation of $I$. Then
\[H^0_m(R/I)=I^{\sat}/I\neq0\]
and so $\depth (R/I)=0$. Suppose that $M(I)<n$. If we consider the
homogeneous ideal $J=I+(x_n)/(x_n)$ in $S=k[x_1,\ldots,x_{n-1}]$
then $M(I)=M(J)$ and
\[
M(J)=(n-1)-\depth (S/J)
\]
by the induction hypothesis. Since $M(I)<n$, $x_n$ is a regular
element of $R/I$ and $\depth (R/I) = 1+\depth (S/J)= n-M(I)$.
\end{proof}

 From Lemma \ref{BR^D(I)}, if we have the generic initial ideal of
a homogeneous ideal $I\subset R$ under a monomial term order
$\tau$ then we may know the dimension of $R/I$. In particular, if
we let $X$ be a closed subscheme and let $I_X$ be the defining
saturated ideal of $X$ in $\mathbb P^{n-1}$, then the codimension
of $X$ is the same as $D(\Gin_{\tau}I_X)$ for any monomial term
order $\tau$. Note that $\Gin_{\tau}(I_X)$ must have a generator
of the form $x_{D(I)}^r$ for some positive number $r$.

\begin{Def}
Let $I$ be a Borel fixed monomial ideal of $R$ and suppose that
\[D(I)<M(I).\]
We denote by $A=\{m_1,\ldots,m_s\}$ the set of minimal monomial
generators of $I$. A monomial $x^K \in R$ is said to be a {\it
generalized sporadic zero} of $I$ if it satisfies the condition that there is a
generator $m_i\in A$ satisfying $\max(m_i)=M(I)$, such that
$x_{M(I)}^rx^K=m_i$ for some $r > 0$. Now we put:
\begin{align*}
&\Spor(I)= \{x^K\in R\,:\, x^K\,\textup{is a generalized sporadic
zero of}\,I\}\\
&\Spor(m,I)=\{x^K\in R\,:\, x^K\in \Spor(I),\, |K|=m\}
\end{align*}
In case $D(I)=M(I)$, we set $\Spor(I)=\emptyset$ for convenience.
\end{Def}

\begin{Coro}
Let $I$ be a Borel fixed monomial ideal. Then the following two
conditions are equivalent:
\begin{enumerate}
\item[(a)] $R/I$ is a graded Cohen-Macaulay ring. \item[(b)]
$D(I)=M(I).$
\end{enumerate}
\end{Coro}

\begin{proof}
This follows immediately from Lemma \ref{BR^D(I)}.
\end{proof}

For a homogeneous ideal $I$, there is a Borel fixed monomial ideal
canonically attached to $I$: the generic initial ideal $\Gin(I)$
with respect to the reverse lexicographic order. It plays a
fundamental role in the investigation of many algebraic,
homological, combinatorial and geometric properties of the ideal
$I$ itself. From now on, we will only use reverse lexicographic
order and we set
\[M(\Gin(I))=M(I),\quad D(\Gin(I))=D(I),\quad \Spor(\Gin(I))=\Spor(I).\]

For many geometric applications, and also for doing inductive
arguments, it is useful to know what happens when we restrict to
a generic hyperplane. Let $h$ be a general linear form. Then we
may consider $J=(I+(h))/(h)$ as a homogeneous ideal of
$S=k[x_1,\ldots,x_{n-1}]$. One has the following well known
fact \cite{BS,G}:
\begin{equation}\label{BS1}
\Gin(J)=\Gin(I)|_{x_n\rightarrow0} .
\end{equation}
Since
$\Gin(I^{sat})=\Gin(I)|_{x_n\rightarrow 1}$ for any homogeneous
ideal $I$, we get
\begin{equation}\label{BS2}
\Gin(J^{sat})=(\Gin(I)|_{\,x_n\rightarrow0})|_{x_{n\!-\!1}\rightarrow
1}.
\end{equation}

\begin{Lem}\label{Sp}
Let $X$ be a closed subscheme in $\bbp^{n-1}$ and let $I_{X}$ be
the defining saturated ideal of $X$. For a general linear form $h$
and for any positive integer $m\geq 1$, we denote by $K(m,I_X)$
the kernel of the map
\[H^1(\mathcal I_X(m-1))\stackrel{\cdot h}{\rightarrow} H^1(\mathcal I_X(m))\]
Then, the following holds:
\begin{equation*}
\dim_kK(m,I_X)=\left\{\begin{array}{ll}
         |\Spor(m,I_X)|&\,\textup{if}\,M(I_X)=n-1;\\
         {}&{}\\
         0&\,\textup{if}\,M(I_X)<n-1.
                 \end{array}\right.
\end{equation*}
\end{Lem}

\begin{proof}
Let $m\geq 1$ and we set
\[
J=\left(\frac{I_{X}+(h)}{(h)}\right).
\]
Then, from (\ref{BS1}), (\ref{BS2}) and from the long exact sequence
\[
0\rightarrow I_X(m-1) \stackrel{\cdot h}{\rightarrow} I_X(m) \rightarrow
I_{X\cap H}(m) \rightarrow H^1(\mathcal I_X (m-1))\stackrel{\cdot
h}{\rightarrow}H^1(\mathcal I_X (m))\rightarrow\cdots
\]
we get the desired result:
\begin{align*}
\dim_k(K(m,I_X))=&\,\dim_k(J^{sat}/J)_m\\
=&\,\dim_k\frac{((\Gin(I_X)|_{\,x_n\rightarrow0})|_{x_{n\!-\!1}\rightarrow
1})_m}{(\Gin(I_X)|_{x_n\rightarrow0})_m}\\
               =&\,\left\{\begin{array}{ll}
         |\Spor(m,I_X)|&\,\textup{if}\,M(I_X)=n-1;\\
         {}&{}\\
         0&\, \textup{otherwise.}
         \end{array}\right.
\end{align*}
\end{proof}

The following theorem gives the relation between the existence of a
generalized sporadic zero and the nonvanishing of sheaf cohomology of
an ideal sheaf $\mathcal I_X$.
\begin{Thm}\label{Cohomology}
Let $X$ be the closed subscheme in $\bbp^{n-1}$ and let $I_{X}$ be
the defining saturated ideal of $X$. Set $d=n-M(I_X)$ and assume
that $\Spor(m,I_X) \neq \emptyset$ for some $m\geq 1$. Then,
\begin{enumerate}
\item[(a)] $H^i(\mathcal I_{X}(j))=0$ for $0<i<d$ and $j\in \mathbb Z.$
\item[(b)] $H^d(\mathcal I_{X}(m-d))\neq 0.$
\end{enumerate}
\end{Thm}

\begin{proof}
To prove this assertion, we use induction on the dimension of the ambient
projective space. Note that $d>0$ if $I_X$ is a saturated ideal.
In case $d=1$, the result follows from Lemma \ref{Sp}.

For $d>1$, we will first prove (a). Note that $M(I_X)<n-1$. By
Lemma \ref{Sp}, the map
\[H^1(\mathcal I_X(j-1))\rightarrow H^1(\mathcal
I_X(j))\] is an inclusion for all integers $j$, and hence
$H^1_\ast(\mathcal I_X)=0$. This gives the proof for $d=2$. Now
assume that $d>2$. From the fact that $M(I_X)=M(I_{X\cap H})$  for a general
hyperplane $H$ (it follows from
Lemma \ref{BR^D(I)}), we have $d-1=(n-1)-M(I_{X\cap H})$. Hence
\[
H^i(\mathcal I_{X\cap H}(j))=0\quad \textup{for}\,0<i<d-1
\]
by the induction hypothesis. Then it follows from the  exact
sequence
  \[
\cdots \rightarrow H^{i-1}(\mathcal I_{X\cap H}(j)) \rightarrow H^i(\mathcal
I_X(j-1))\rightarrow H^i(\mathcal I_X(j))\rightarrow
\]
  that $H^i(\mathcal I_{X}(j))=0$ for $0<i<d$ and for all $j$. This
  proves (a).

We now prove (b). For $d>1$, note again that $M(I_X)=M(I_{X\cap H})$.
Hence we have
\[
H^{d-1}(\mathcal I_{X\cap H}(m-d+1))\neq 0
\]
by the induction hypothesis. From the  exact sequence
\[
0\rightarrow H^{d-1}(\mathcal I_{X\cap H}(m-d+1))
\rightarrow H^d(\mathcal I_X(m-d))\rightarrow H^d(\mathcal
I_X(m-d+1))
\]
we know that $H^d\big(\mathcal I_X(m-d)\big)$ does not
vanish, and this proves (b).
\end{proof}

 From this, we give a cohomological proof of the following well
known fact.
\begin{Coro}\textup{(\cite{E1})}
For $I$ be a homogeneous ideal of $R=\kxn$, consider the generic
initial ideal for the degree reverse lexicographic order. Then,
\[
\depth (R/I)=\depth \left(R/\Gin(I)\right).
\]
\end{Coro}

\begin{proof}
Let $H^i_m(R/I)$ be the local cohomology of $R/I$. By Lemma
\ref{BR^D(I)}, we know that
\[
\depth
\left(R/\Gin(I)\right)=n-M(I),
\]
so it is enough to show that
\[
n-M(I)=\min\{i:\,H^i_m(R/I)\neq 0\}.
\]
But this follows from Theorem \ref{Cohomology} and the following
fact (\cite{E2}):
\[
H^i_m(R/I) = \bigoplus_{j \in \mathbb Z }H^i(\bbp^n,\mathcal I(j))
\]
for $0<i<n$.
\end{proof}

\begin{Coro}\label{Cohen}
Let $I$ be a homogeneous ideal of $R=\kxn$. Then the following
facts are equivalent.
\begin{enumerate}
\item[(a)] $R/I$ is Cohen-Macaulay.
\item[(b)] $R/\Gin(I)$ is Cohen Macaulay.
\item[(c)]$M(I)=D(I)$.
\end{enumerate}
\end{Coro}

\begin{Ex}
We consider the homogeneous ideal in $k[x_1,x_2,x_3,x_4]$
     \[I=(x_3^3-x_1x_4^2, x_1^2x_3^2-x_2^3x_4, x_2^3x_3 - x_1^3x_4,
x_2^6-x_1^5x_3).\]
Using Macaulay2, we get the generic initial ideal of $I$
\[\Gin(I)=(x_1^3, x_1^2x_2^2, x_1x_2^3, x_2^5, x_2^4x_3^2) \]
under reverse lexicographic order. Then we know that $I$ is a
saturated ideal, $D(I)=2$ and $M(I)=3$. Hence $I$ is the defining
ideal of a projective curve $C$ in $\mathbb P^3$. Note that
$\Spor(I)=\{x_2^4, x_2^4x_3\}$ and $C$ is not arithmetically
Cohen-Macaulay.
\end{Ex}

\begin{Remk}
  If we use another monomial term order then $M(I)$ may change, but
  $D(I)$ is independent of the monomial order.
\end{Remk}

By Theorem \ref{BS}, the regularity of $\Gin(I)$ is the largest
degree of a generator of $\Gin(I)$. Bayer and Stillman \cite{BS}
showed the regularity of $I$ is equal to the regularity of
$\Gin(I)$.
\begin{Thm}\label{BS_reg}
{\rm (\cite{BS}, \cite{G}, \cite{GS})} For any homogeneous ideal
$I$, using the reverse lexicographic order,
\[\sat(I)=\sat(\Gin(I)),\,\]
\[\reg(I)=\reg(\Gin(I)).\]
\end{Thm}

For a homogeneous ideal $I\subset R$ there exists a flat family of
ideals $I_t$ with $I_0=\ini(I)$ and $I_t$ canonically isomorphic
to $I$ for all $t\neq 0$ (Corollary 1.21 in \cite{G}). Using this
result we know the minimal free resolution of $I$ is obtained from
that of $\ini(I)$ by cancelling some adjacent terms of the same
degree. That is, we can always choose a complex of $K\cong
R/m-$modules $V_{\bullet}^d$ such that
\[V_i^d\cong \Tor_{i}^R(\ini(I),K)_d\]
\[H_i(V_{\bullet}^d)\cong \Tor_{i}^R(I,K)_d.\]

Although the minimal free resolution of a general monomial ideal
can be quite complicated, the situation for Borel fixed monomial
ideals is very nice. Eliahou and Kervaire \cite{EK} gave a result
for the structure of the minimal free resolution of a Borel fixed
ideal. Using these results, we have the following.
\begin{Thm}\label{CP}
\textup{\bf(Crystallization Principle)} Let $I$ be a homogeneous
ideal generated in degrees $\leq d$. Assume that there is a
monomial order $\tau$ such that $\Gin_{\tau}(I)$ has no generator
in degree $d+1$.  Then $\Gin_{\tau}(I)$ is generated in degrees
$\leq d$ and $I$ is $d$-regular.
\end{Thm}

\begin{proof}
The case of arbitrary monomial order $\tau$ can be proved in the
same manner as the proof of Theorem 2.28 in \cite{G}.
\end{proof}

\begin{Remk}
We can consider Theorem \ref{CP} as a generalization of the
Gotzmann persistence theorem (Theorem 3.8 in \cite{G}). Let
$I^{\textup{lex}}$ be the lex-segment ideal of $I$. We know that
the lex-segment ideal has the largest Betti numbers in the class
of the ideals with a given Hilbert function (Theorem 2 in
\cite{B}). Then $\Gin_{\tau}(I)$ has no generator in degree $d+1$
for every monomial order $\tau$ if $I^{\textup{lex}}$ has no
generator in degree $d+1$, which is equivalent to the Hilbert
function of $R/I$ has maximal growth in degree $d$. Hence if a
homogeneous ideal $I$ is generated in degrees $\leq d$ and the
Hilbert function of $R/I$ has maximal growth in degree $d$ then
$I$ is $d$-regular. We will use Theorem \ref{CP} to generalize
results in \cite{BGM}.
\end{Remk}

For a homogeneous ideal $I\subset R=\kxn$, let $d=\dim(R/I)$. L.T.
Hoa and N. V. Trung defined the {\it s-reduction number} of $R/I$
for $s\geq d$ in \cite{HT}. They have shown that
$r_s(R/I)=r_s(R/\Gin(I))$ (Theorem 2.3 in \cite{HT}). If $I$ is a
Borel fixed monomial ideal we know that there are positive numbers
$t_1,\ldots,t_d$ such that $x_i^{t_i}$ is a minimal generator of
$\Gin(I)$. In \cite{HT}, authors also proved that if a monomial
ideal $I$ is strongly stable, then
\begin{equation}\label{eq4}
r_s(R/I)=\min\{k\,:\,x_{n-s}^{k+1}\in I\}.
\end{equation}
 From these facts we get the following Lemma.  We will take this lemma
as the definition of $r_s(R/I)$ for our purposes.

\begin{Lem}\label{Reduction Number}
For a homogeneous ideal $I$ of $R$ and for $s\geq \dim(R/I)$,
the $s$-reduction number $r_s(R/I)$ can be given as the following:
\begin{align*}
r_s(R/I)=&\min\{k\,:\,x_{n-s}^{k+1}\in \Gin(I)\}\\
         =&\min\{k\,: \textup{
Hilbert function of }R/(I+J) \textup{ vanishes in degree }k+1 \textup{ }\}
\end{align*}
where $J$ is generated by $s$ general linear forms of $R$.
\end{Lem}

\begin{proof}
The first part of Lemma follows directly from results in \cite{HT}.
Consider a homogeneous ideal $J$ generated by $s$ general linear
forms. Then, for the reverse lexicographic order, we know that
\[M:=\Gin\left(\frac{I+J}{J}\right)=\frac{\Gin(I)+(x_{n-s+1},\ldots,x_n)}{(x_{n-s+1},\ldots,x_n)}\]
from (\ref{BS1}). Hence, for $s\geq d$, we can compute
$r_s({R/I})$ from (\ref{eq4}) by looking for the smallest degree
in which the monomial ideal $M$ becomes the unique maximal ideal of
$R/(x_{n-s+1},\ldots,x_n)$ by the strongly stable property of $M$.
This proves the second part of the Lemma since the Hilbert function of
$(I+J)/J$ has the same Hilbert function as $M$.
\end{proof}

\begin{Remk} \label{WLP connection}
Let $\Gamma$ be a finite set of points (or more generally, any
zero- \linebreak dimensional scheme).  There is a strong connection between the
definition of \linebreak
$r_2(R/I_\Gamma)$ and the Weak Lefschetz Property (WLP) (cf.\ \cite{HMNW}). We
first recall this property.  Let $J = (L_1,L_2)$ be generated by general
linear forms.  Let $K = I_\Gamma + (L_1)$, and let $A = R/K$ be the Artinian
reduction of $R/I_\Gamma$ by $L_1$.  Then multiplication by $L_2$
gives an exact
sequence
\[
0 \rightarrow
\left ( \frac{[K:L_2]}{K} \right )_d \rightarrow
(R/K)_d \stackrel{\times L_2}{\longrightarrow}
(R/K)_{d+1} \rightarrow R/(I_\Gamma + J)_{d+1} \rightarrow 0.
\]
The Weak Lefschetz Property merely says that this multiplication has maximal
rank, for all $d$.  Now, computing $r_2(R/I_\Gamma)$ amounts to studying the
surjectivity of $\times L_2$ above, which is a triviality when we know maximal
rank.  Indeed, WLP says that there is an expected value for
$r_2(R/I_\Gamma)$, and
that this value is achieved.

Intuitively, WLP should be true ``most of the time."  Attempts to
make this more
precise were given in \cite{HMNW} and in \cite{MMR3}.  A very natural question,
then, is whether WLP should hold for points with the Uniform Position Property
(UPP).  In  Example \ref{UPP vs WLP} we will show that this is not the case.
\end{Remk}

In the case of low dimensional schemes, the following Lemma
\ref{regularity_reduction} and Proposition \ref{firstCH regular}
show a relation between the reduction number and the regularity.

\begin{Lem}\label{regularity_reduction}
Let $I\subset R=\kxn$ be a saturated ideal and suppose that
$\dim(R/I)=1$. Then the regularity of $I$ is equal to
$r=r_1(R/I)+1$.
\end{Lem}

\begin{proof}
Consider the generic initial ideal $\Gin(I)$ in terms of reverse
lexicographic order. Then $x_{n-1}^{r}$ is a minimal generator of
$\Gin(I)$ (Lemma \ref{Reduction Number}). We know that all
monomials of degree $r$ that do not involve the variable $x_n$ are contained in
$\Gin(I)$ by the strongly stable property. Hence if there is a
minimal generator of
$\Gin(I)$ in degree $> r$ then it should involve the variable $x_n$.
But this is
impossible because $\Gin(I)$ is a saturated Borel fixed monomial ideal (Theorem
\ref{BS} and Theorem \ref{BS_reg}). This implies that the maximal
degree of a minimal generator of $\Gin(I)$ is $r$, which is equal
to the regularity of $I$ by Theorem \ref{BS_reg}.
\end{proof}

\begin{Prop}\label{firstCH regular}
Let $C$ be a projective curve in $\mathbb P^{n-1}$, $n>3$ (not
necessarily reduced and irreducible). Suppose that there is an integer $d>
r_2(R/I_{C})$ such that
\[h^1(\mathcal I_C(d-1))=0.\]
Then $I_C$ is $d$-regular.
\end{Prop}

Proposition \ref{firstCH regular} shows when we can get information about the
regularity of a projective curve $C\subset \mathbb P^{n-1}$ by the
vanishing of $h^1(\mathcal I_C(t))$.  Normally it requires some knowledge of
$h^2(\mathcal I_C(t))$ as well.

\begin{proof}
We will prove that $h^2(\mathcal I_C(d-2))=0$. For a general
hyperplane $H$,
\[r_1(I_{C\cap H})\leq r_2(I_C)\]
since $\Gin(I_{C\cap H})=(\Gin(I_C)|_{x_n\rightarrow
0})|_{x_{n-1}\rightarrow 1}$ (see Proposition \ref{Reduction
Number}).
  Therefore we have that
  \[\reg (I_{C\cap H}) \leq
r_2(R/I_{C})+1\]
by Lemma \ref{regularity_reduction}, and so
\[h^1(\mathcal I_{C\cap H}(t))=0\]
for all $t\geq r_2(R/I_{C})$. Consider the following long exact
sequence
\[\longrightarrow(h^1(\mathcal I_{C\cap
H}(t))\longrightarrow h^2(\mathcal I_C(t-1))\longrightarrow
h^2(\mathcal I_C(t))\longrightarrow, \] then we get
\[0\longrightarrow h^2(\mathcal I_C(t-1))\longrightarrow h^2(\mathcal I_C(t))\]
is injective for all $t\geq r_2(R/I_{C})$ and so $h^2(\mathcal
I_C(t-1))=0$ for all $t\geq r_2(R/I_{C})$. This implies that
$h^2(\mathcal I_C(d-2))=0$ and $I_C$ is $d$-regular.
\end{proof}


\section{The regularity of a Cohen-Macaulay ring} \label{reg sect}
Let $X$ be a closed subscheme of $\bbp^{n-1}$ which is not
contained in any hyperplane. Suppose that $I_{X}$ is the defining
saturated ideal of $X$ and that $R/I_X$ is Cohen-Macaulay. Then
the following is a well-known fact (\cite{E2} Corollary 4.14, \cite{EG}):
\[
\reg(I_X)\leq \deg(X)-\codim(X)+1.
\]
In this section, we will give a stronger bound,
by using  the initial degree of $I_X$ (Theorem \ref{Cohen1}). We start with
the following:

\begin{Prop}\label{Re}
Let $I$ be a saturated Borel fixed monomial ideal of $R=\kxn$ with
$\dim (R/I)=1$. Then,
\[
\reg \ I \leq \deg (R/I) - \binom{\alpha-1+n-1}{\alpha-1}+\alpha
\]
where $\alpha= \min \{l \mid I_l \neq 0 \}$ is the initial degree
of $I$.
\end{Prop}

\begin{proof}
Since $R/I$ is a Cohen-Macaulay ring,
\[
D(I)=M(I)=n-1.
\]
Hence, there exist positive numbers $\lambda_1\leq
\lambda_2\leq\cdots\leq\lambda_{n-1}$ such that
\[
I=(x_1^\alpha,\ldots,x_2^{\lambda_1},\ldots,x_{n-2}^{\lambda_{n-2}},
\ldots,x_{n-1}^{\lambda_{n-1}}).
\]
Clearly if $J=(x_1,\ldots,x_{n-1})^\alpha$, then
\[
\deg R/J
=\binom{\alpha+n-1}{\alpha}-\binom{\alpha+n-2}{\alpha}=
\binom{\alpha-1+n-1}{\alpha-1}.
\]
Now, let $J'$ be the ideal generated by
$x_{n-1}^{\lambda_{n-1}}$ and by all the monomials of degree
$\alpha$ except $ x_{n-1}^\alpha$ in the variables
$x_1\ldots,x_{n-1}$:
\[
J'=(x_1^\alpha,\ x_1^{\alpha-1}x_2,\ \ldots,x_{n-2}x_{n-1}^{\alpha-1},\
x_{n-1}^{\lambda_{n-1}}).
\]
We first claim that
\begin{equation}\label{Re1}
\deg (R/J')=\binom{\alpha-1+n-1}{\alpha-1}-\alpha+ \lambda_{n-1}.
\end{equation}
Indeed, if we consider the first difference Hilbert function (i.e
$h$-vector) of $R/J'$, then $\Delta H(R/J', t)=1$ if $\alpha \leq
t < \lambda_{n-1}$, since $J'$ is defined by using all the
monomials of degree $\alpha$ except $x_{n-1}^{\alpha}$ and we have
reduced modulo a general linear form (which we can take to be
$x_n$). Since
\[
H(R/J',\alpha)=H(R/J,\alpha)+1=\binom{\alpha-1+n-1}{\alpha-1}+1,
\]
(\ref{Re1}) follows from
\[
\deg(R/J')=H(R/J',\lambda_{n-1}-1)=H(R/J',
\alpha)+\sum_{t=\alpha+1}^{\lambda_{n-1}-1}(\Delta H(R/J', t)).
\]

Note also that
\begin{equation}\label{Re2}
\deg (R/J') \leq \deg (R/I).
\end{equation}
Since $I$ and $J'$ are also Borel fixed monomial ideals, by virtue
of Theorem \ref{BS},
\[
\lambda_{n-1}=\reg (I)=\reg(J')
\]
and it follows from (\ref{Re1}) and (\ref{Re2}) that
\[
\reg(I)\leq \deg (R/I)-\binom{\alpha-1+n-1}{\alpha-1}+\alpha
\]
as claimed.
\end{proof}

\begin{Coro}\label{Re3}
If $I$ is a Borel fixed monomial ideal such that $R/I$ is
Cohen-Macaulay with $\dim R/I \geq 1$ then
\[\reg(I)  \leq \deg(R/I)-\binom{\alpha-1+D(I)}{\alpha-1}+\alpha\]
where $\alpha=\alpha(I)= \min \{l \mid I_l \neq 0 \}$.
\end{Coro}

\begin{proof}
To get the regularity bound, we work by induction on $n$. If
$D(I)=n-1$ then it is proved by Proposition \ref{Re}. Now suppose
that $D(I)<n-1$. If we consider the ideal $J=I+(x_n)/(x_n)$ of
$S=k[x_1,\ldots,x_{n-1}]$, then $J$ is also a Borel fixed monomial
ideal such that
\[\reg(I) =\reg(J) ,\quad D(I)=D(J)\]
and
\[\deg(R/I)= \deg(S/J),\quad \alpha(I)=\alpha(J).\]
Since $S/J$ is also a Cohen-Macaulay ring, by the induction
hypothesis,
\begin{align*}
\reg(I) = \reg(J) &\leq
\deg(R/J)-\binom{\alpha(J)-1+D(J)}{\alpha(J)-1}+\alpha(J)\\
                 &= \deg(R/I)-\binom{\alpha(I)-1+D(I)}{\alpha(I)-1}+\alpha(I)
\end{align*}
as desired.
\end{proof}

The following is the main result of this section.  It is a special case of Satz
3 of \cite{uwe}, taking $I(A) = 0$ since $X$ is arithmetically Cohen-Macaulay. 
However, our method of proof is quite different.

\begin{Thm}\label{Cohen1}
Let $X$ be an arithmetically Cohen-Macaulay subscheme of $\bbp^{n-1}$  with
$\codim (X)=e$. Then,
\[\reg (I_X)\leq\deg (X) -\binom{\alpha-1+e}{\alpha-1}+\alpha\]
where $\alpha=\alpha(I_X)$ is the initial degree of $I_X$.
\end{Thm}

\begin{proof}
It follows from Lemma \ref{BR^D(I)}, Corollary \ref{Cohen} and
Corollary \ref{Re3}, since $\alpha(I_X)=\alpha(\Gin(I_X))$.
\end{proof}

\begin{Remk}
Under the assumption of Theorem \ref{Cohen1}, consider the
function $f(x) = $ \linebreak $\binom{x-1+e}{x-1}-x$. If there is
no linear form in $I_X$, we have:
\begin{align*}
   \reg(I_X)&\leq \deg(R/I_X)-f(\alpha(I_X))\\
          &\leq \deg(R/I_X)-f(2)\\
          &=\deg(R/I_X)-D(I_X)+1\\
          &=\deg(R/I_X)-\codim(R/I_X)+1.
\end{align*}
since $f(x)$ is strictly increasing for $x \geq 2$. Hence we may
think of Theorem \ref{Cohen1} as a generalization of a result in
\cite{E2}.
\end{Remk}


\section{ Application to the Hilbert functions of a set of points} \label{first
points sect}

  A.\ Bigatti, A.V.\ Geramita and J.C.\ Migliore gave many
geometric results relating to the maximal growth of the first difference
of Hilbert function in \cite{BGM}. In particular, they considered a
reduced set of points $\Gamma$ in $\mathbb P^{n-1}$ and proved
several geometric
consequences of the condition
\begin{equation}\label{eq2}
\Delta H(\Gamma,d)=\Delta H(\Gamma,d+1)=s
\end{equation}
for $d\geq s$. In this
section, we generalize these results of \cite{BGM} by assuming only $d >
r_2(R/I_Z)$. In the next section we will take up the question of
uniform position,
and see what can be said (and what can {\em not} be said) with that additional
hypothesis.  We will see that before we add the uniformity
assumption, virtually
everything that holds in the case $d \geq s$ (\cite{BGM} Theorem 3.6)
remains true
here (Theorem \ref{Generalize}).  The only difference will be our
inability here
to prove a reducedness result (see Remark \ref{compare BGM}).  We
start with the
following.

\begin{Prop}\label{Decreasing}
Let $Y$ be a scheme of dimension $\leq 1$ in $\mathbb P^{n-1}$,
$n\geq 3$. Then
\begin{equation}\label{eq1}
\Delta H(R/I_{Y},d)\geq \Delta H(R/I_{Y},d+1)
\end{equation}
for $d\geq r_2(R/I_{Y})$. Moreover, for $d> r_2(R/I_{Y})$, we
have equality in {\rm({\ref{eq1}})} if and only if $\Gin (I_Y)$
has no minimal generators in degree $d+1$.
\end{Prop}

\begin{proof}
Let $J = (L_1,L_2)$ be generated by general linear forms. Let $K =
I_Y + (L_1)$; this can be viewed as a homogeneous ideal in $S=R/(L_1)$.  Let $A =
S/K$. Then
multiplication by $L_2$ gives an exact sequence
\begin{equation} \label{eq8}
0 \rightarrow \left ( \frac{[K:L_2]}{K} \right )_d \rightarrow
(S/K)_d \stackrel{\times L_2}{\longrightarrow} (S/K)_{d+1}
\rightarrow (R/(I_\Gamma + J))_{d+1} \rightarrow 0.
\end{equation}
By the definition of $r_2(R/I_{Y})$,
\[
\left (R/(I_\Gamma + J) \right )_{d+1}=0
\]
for $d\geq r_2(R/I_Y)$. Since $I_Y$ is a saturated ideal, $\Delta
H(R/I_{Y},d)$ is the same as the Hilbert function $H(S/K,d)$ for
all $d$. This proves the first part of the theorem.

Now consider the following exact sequence:
\begin{equation}\label{exact_sequence}
0 \rightarrow \left ( \frac{[K:L_2]}{K} \right )_d \rightarrow
(S/K)_d \stackrel{\times L_2}{\longrightarrow}
(S/K)_{d+1}\rightarrow 0
\end{equation}
for $d> r_2(R/I_{Y})$. Note that $r_s(R/I_{Y})=r_s(R/\Gin(I_{Y}))$
by Theorem 2.3 in \cite{HT}. So we may replace $K$, $L_1$ and
$L_2$ by $\Gin(K)$, $x_{n}$ and $x_{n-1}$ respectively and thus
reduce to the case where $K$ is a Borel fixed monomial ideal,
since we know that (\ref{BS1}) and (\ref{BS2}) hold under reverse
lexicographic order (note that $K$ and $\Gin(K)$ have the same
Hilbert function).

Since $I_{Y}$ is a saturated ideal, $\Gin(K)$ have no minimal generator of
degree $d+1$ in $S=R/L_1$ if and only if $\Gin(I_Y)$ have no minimal generator
of degree $d+1$ in $R$ for $d>r_2(R/I_Y)$.  For the proof of
the second part, therefore, it is enough to show that $\Gin (K)$
has no minimal generators in degree $d+1$ if and only if
\begin{equation}\label{eq5} \left (
\frac{[\Gin(K):x_{n-1}]}{\Gin(K)} \right )_d=0
\end{equation}
from the exact sequence (\ref{exact_sequence}).

Now assume that (\ref{eq5}) holds. We know that $x_{n-2}^{d}\in \Gin(K)$
by the first equality of Lemma \ref{Reduction Number}. This
implies that $\Gin(K)_d$ contains each monomial $x^{J}$ of degree
$d$ such that
\[
\supp(J)\subset \{x_1,\ldots,x_{n-2}\}
\] 
since $\Gin(K)$ is a Borel fixed
monomial ideal. If $\Gin(K)$ has a minimal generator of degree
$d+1$, then it must involve the variable $x_{n-1}$. But this is
impossible because it means $\left ( {[\Gin(K):x_{n-1}]}/{\Gin(K)}
\right )_d\neq 0$. Hence $\Gin (K)$ has no minimal generators in
degree $d+1$. 

Conversely, if
\[
\left
( {[\Gin(K):x_{n-1}]}/{\Gin(K)} \right )_d\neq 0
\] 
then there
exists a monomial $x^L$ such that $x_{n-1}x^L \in \Gin(K)_{d+1}$
but $x^L \notin \Gin(K)_{d}$. Assume that $x_{n-1}x^L$ is not a
minimal generator of $\Gin(K)$. Then we can choose a monomial
$x^J\in \Gin(K)_d$ of degree $d$, satisfying
\[x_ix^J=x_{n-1}x^L\]
for $1 \leq i<{n-1}$. Note that $x^J$ should contain the variable
$x_{n-1}$. Hence we have
\[x^L=x^J\left(\frac{x_i}{x_{n-1}}\right) \in \Gin(K)_{d}\]
by the Borel fixed property (or strongly stable) and this contradicts the choice
of $x^L$. Hence the monomial $x_{n-1}x^L$ is a minimal generator of degree
$d+1$ of $\Gin(K)$.
\end{proof}

\begin{Thm}\label{Saturated1}
Let $Y$ be a scheme of dimension $\leq 1$ in $\mathbb P^{n-1}$,
$n\geq 3$. Assume that for some $d>r_2(R/I_Y)$
\[
\Delta H(R/I_{Y},d)= \Delta H(R/I_{Y},d+1))=s\neq 0.
\]
Then $\langle (I_Y)_{\leq d}\rangle$ is $d$-regular and it is a
saturated ideal defining a one dimensional scheme of degree $s$ in
$\mathbb P^{n-1}$.
\end{Thm}

\begin{proof}
Let $\bar I=\langle (I_Y)_{\leq d}\rangle$. To show that $\bar I$
is a saturated ideal, we need only to prove that $\Gin(\bar I)$
does not have a minimal generator containing the variable $x_{n}$
(Theorem \ref{BS}). By Theorem \ref{CP} and Theorem
\ref{Decreasing}, $\bar I$ is $d$-regular, and thus $\Gin(\bar I)$
has no minimal generator in degrees $>d$. Since $I_Y$ is a saturated
ideal, $\Gin(I_Y)$ does not have a minimal generator involving
$x_{n}$. Then it follows from
\[(\Gin(I_Y))_k=(\Gin (\bar I))_k\] for $k\leq d$ that there is no minimal
generator
of $\Gin (\bar I)$ containing the variable $x_n$ in degree $\leq d$.
Hence $\bar I$ is a saturated ideal.

Since $(\bar I)_i=(I_Y)_i$ for $i\leq d+1$, we know that
\[r_2(R/\bar I)=r_2(R/I_Y)\]
and that $\bar I$ defines a scheme of dimension $\leq 1$. Note
that $\Gin(\bar I)$ does not have minimal generator of degree
$k>d$. Applying Theorem \ref{Decreasing} repeatedly, we get
\[\Delta H(R/\bar I,i)= \Delta H(R/\bar I,i+1))=s\neq 0\]
for all $i\geq d> r_2(R/\bar I)$. This means that $\bar I$ defines a
one dimensional subscheme in $\mathbb P^{n-1}$ of degree $s$.
\end{proof}

\begin{Coro}\label{Saturated}
Let $\Gamma$ be a set of points in $\mathbb P^{n-1}$, $n\geq 3$,
and let $(1,h_1,\ldots,h_t)$ be the $h$-vector of $R/I_{\Gamma}$.
Then $h_d\geq h_{d+1}$ for $d>r_2(R/I_{\Gamma})$. Suppose that
\[h_{d}=h_{d+1}=s\] for $d> r_2(R/I_{\Gamma})$. Then
$(I_{\Gamma})_{\leq d}$ is a saturated ideal defining a one
dimensional subscheme of degree $s$ in $\mathbb P^{n-1}$ and it is
$d$-regular.
\end{Coro}

\begin{proof}
Since $\Gamma$ is arithmetically Cohen-Macaulay,
\[\Delta H(R/I_{\Gamma},d)=h_d.\]
Hence the result follows from Corollary \ref{Saturated1}.
\end{proof}

\begin{Remk}\label{remark1}
Let $Y\subset \mathbb P^{n-1}$ be a subscheme of any dimension,
$n\geq 3$. Assume that
\[\Delta H(Y,d)=\Delta H(Y,d+1)=s\neq 0\]
for some $d\geq s$. Then $d>r_2(R/I_{Y})$ and $\dim(Y)\leq 1$.
Indeed, the condition $d \geq s$ means that the $d$-binomial expansion
of $s$ is $\binom{d}{d}+\cdots+\binom{d-s+1}{d-s+1}$ and so we
have maximal growth of the Hilbert function of $R/(I_{Y}+(L))$ in
degree $d$, where $L$ is a general linear form of $R$. By Theorem 3.6
of \cite{BGM},
this implies that $(I_{Y})_{\leq d}$ is a saturated ideal defining a
one dimensional
subscheme in $\mathbb P^{n-1}$. Then
\[D((I_{Y})_{\leq d})=n-2\]
by Lemma \ref{BR^D(I)}, and $x_{n-2}^{d}$ is contained in
$(I_{Y})_{\leq d}$. This implies that $d> r_2(R/I_{Y})$ and that
$D(I_Y)\geq n-2$, so that $\dim(Y)\leq 1$. Hence we can think of the
condition $d>r_2(R/I_{Y})$ as generalizing the
assumption $d \geq s$ which is in \cite{BGM}.
\end{Remk}

\begin{Remk}
  \label{not curve}
Since Corollary \ref{Saturated} weakens the hypothesis $d \geq s$ of
\cite{BGM},
one might ask if our hypothesis $d > r_2(R/I)$ can be weakened even
further.  We now
show by simple examples that this is not the case for results involving the
first difference, although in section \ref{add'l comments} we do weaken this
hypothesis by invoking the second difference.

Let $\Gamma$ be a complete intersection in $\mathbb P^3$ of type
$(4,4,4)$.  The
first difference of the Hilbert function of $\Gamma$ is
\[
1 \ \ 3 \ \ 6 \ \ 10 \ \ 12 \ \ 12 \ \ 10 \ \ 6 \ \ 3 \ \ 1.
\]
Since the Artinian reduction of $R/I_\Gamma$ has the Weak Lefschetz property
(\cite{HMNW} Theorem 2.3), we see that $r_2(R/I_\Gamma) = 4$.  Since $\Delta
H(R/I_\Gamma,4) = \Delta H(R/I_\Gamma,5) = 12$,
but clearly the component $(I_\Gamma)_4$ defines a zeroscheme (namely $\Gamma$)
rather than a curve, we see that we cannot weaken the assumption $d >
r_2(R/I_\Gamma)$.  Of course such an example can be constructed in any
codimension $\geq 3$.

Similarly, we have the following series of examples (and again, the
same sort of
thing could be done in higher projective space):
\[
\begin{array}{l}
\hbox{7  general points in $\mathbb P^3$ have $h$-vector 1 \ 3 \ 3} \\
\hbox{16 general points in $\mathbb P^3$ have  $h$-vector 1 \ 3 \ 6 \ 6} \\
\hbox{30 general points in $\mathbb P^3$ have  $h$-vector 1 \ 3 \ 6 \
10 \ 10} \\
\hbox{etc.}
\end{array}
\]
In each case, the value of $r_2(R/I_\Gamma)$ is the ``expected'' one,
and we have $d
= r_2(R/I_\Gamma)$, but the component of $I_\Gamma$ in degree $d$ is zero,
while the component of $I_\Gamma$ in degree $d+1$ defines a zero-dimensional
scheme rather than a curve.  See also Example
\ref{UPP vs WLP}.

We also remark that in general
$r_2(R/I_Z)$ can be much smaller than $s$.  For example, let $Z$ be a set of
sufficiently  many points on a smooth arithmetically Cohen-Macaulay
curve of degree
28 in ${\mathbb P}^3$ with $h$-vector $(1,2,3,4,5,6,7)$.  Then it is
not hard to
show that $r_2(R/I_Z) = 6$ while $s = 28$.
\end{Remk}

The next result is our analog to \cite{BGM} Corollary 3.7.

\begin{Thm}\label{Generalize}
Let $Y\subset \mathbb P^{n-1}$, $n\geq 3$ be a reduced scheme of
any dimension. Assume that $\Delta H(Y,d)=s$ and that the
saturated ideal $\langle (I_Y)_{\leq d}\rangle^{\sat}$ defines a
curve $V$ of degree $s$. Then
\begin{enumerate}
\item[(a)] $\dim(Y)\leq 1$. \item[(b)] $\langle (I_Y)_{\leq
d}\rangle$ is the saturated ideal of a curve $V$ of degree $s$
(not necessarily unmixed).
\end{enumerate}

\noindent Let $C$ be the unmixed one-dimensional component of the curve $V$
in $(b)$. Let $Y_1$ be the subvariety of $Y$ on $C$ and $Y_2$ the
``residual" subvariety.
\begin{enumerate}
\item[(c)]
$\langle (I_{Y_1})_{\leq d}\rangle=I_C$.
\item[(d)] $\dim(Y_2)=0$
and $H(Y_1,t)=H(Y,t)-|Y_2|$ for all $t\geq d-1$.
\item[(e)]
\[
\Delta H(Y_1,t)=\left\{
                                    \begin{array}{ll}
                                    \Delta H(C,t),&\,\,\mathrm{for}\, t\leq
d+1\\
                                    \Delta H(Y,t),&\,\,\mathrm{for}\,
                                    t\geq d.
                                    \end{array}\right..
\]
In particular, $\Delta
H(Y_1,t)\geq \Delta H(Y_1,t+1)$ for all
$t \geq d$. \item[(f)] If we assume $h^1(\mathcal
I_{C_\red}(d-1))=0$ then $V$ is reduced and $C=C_{red}$ is
$d$-regular.
\end{enumerate}
\end{Thm}

\begin{proof}
The proof is almost same as that of Theorem 3.6 in \cite{BGM}. We
  sketch the proof (see \cite {BGM} for more details). Let $\bar I=\langle
(I_{Y})_{\leq d}\rangle$. If $L$ is a general linear form and $H$
is the corresponding hyperplane, let
\[J=\bar I+(L)/(L).\] Then $J^{\sat}$ is the defining ideal of the set of
points $I_{V\cap H}$ of degree $s$ in $S=R/(L)$. By the assumption that
$\Delta H(Y,d)=s$, we get
\[\Delta H(R/\bar I,d)=H(S/J,d)=\deg (S/J)=s.\] This implies $d\geq \reg(J)$.
Hence $J$ is
saturated in degree $\geq d$ and $H(S/J,k)=s$ for all $k\geq d$.
Then, by Lemma \ref{regularity_reduction}
\[d\geq \reg(J)\geq \reg(J^{\sat})
=r_1(S/J^{\sat})+1\] and thus $x_{n-2}^{d}$ is contained in
$\Gin(J^{\sat})_d=\Gin(J)_d$. Hence $d\geq r_1(S/J)+1=r_2(R/\bar
I)+1$ so $\bar I=\langle (I_{Y})_{\leq d}\rangle$ is a saturated
ideal defining a one dimensional scheme $V$ of degree $s$ in
$\mathbb P^{n-1}$ by Corollary \ref{Saturated1}. Then $\dim(Y)\leq
1$ follows from $D(I_Y)\geq n-2$. This proves (a) and (b). Let $C$
be the unmixed one-dimensional component of $V$ with degree $s$.
Then the scheme $Y$ is the union of $Y_1$, which is the subvariety
of $Y$ lying on $C$, and $Y_2$, the remaining subvariety of $Y$.
Since $I_V=\langle (I_{Y})_{\leq d}\rangle$ is $d$-regular
(Theorem \ref{Saturated1}), the same proof as in Theorem 3.3 and
Theorem 3.6 of \cite{BGM} gives
\begin{align}
\label{blah1}   h^1(\mathcal I_{Y_2}(t))&=0 \textup{ for all } t\geq d-1\\
\label{blah2} h^1(\mathcal I_{C\cup Y_2}(t))&=0 \textup{ for all } t\geq d-1\\
\label{blah3} h^1(\mathcal I_{C}(t))&=0 \textup{ for all } t\geq d-1.
\end{align}
Since $Y_2$ is a zero dimensional scheme, we have that $\Delta
H(Y_2,t)=0$ for all $t\geq d$. We also have
\begin{equation}\label{eq3}
H(C\cup Y_2,t)=H(Y,t)+H(C,t)-H(Y_1,t)
\end{equation}
from the short exact sequence
\[
0\longrightarrow(I_{C\cup
Y_2})_t\longrightarrow(I_Y)_t\oplus(I_C)_t\longrightarrow(I_{Y_1})_t\longrightarrow
0.
\]
On the other hand, from the exact sequence
\[
0\longrightarrow\mathcal I_{C\cup Y_2}\longrightarrow\mathcal O_{\mathbb
P^{n-1}}\longrightarrow\mathcal O_{C\cup Y_2}\longrightarrow 0
\]
(respectively the same sequence with $C\cup Y_2$ replaced by $C$),
we get from (\ref{blah2}) and (\ref{blah3}) that $H(C\cup Y_2,t)=h^0(\mathcal
O_C(t))+|Y_2|$ and $h^0(\mathcal O_C(t))=H(C,t)$ for $t\geq d-1$.
Now combining this information and substituting it in (\ref{eq3})
we get
\[
H(Y_1,t)=H(Y,t)-|Y_2|
\]
for all $t\geq d-1$. This proves (d). In particular,
\[
\Delta H(Y_1,t)=\Delta H(Y,t) \textup{ for all }t\geq d.
\]
And this proves the second half of the Hilbert function claimed in
(e). Notice that we still have the following condition for $Y_1$:
\[
\Delta H(Y_1,d)=\Delta H(Y_1,d+1)
\]
Since $I_{Y}\subset I_{Y_1}$, clearly $r_2(R/I_{Y_1})\leq
r_2(R/I_{Y})<d$. Hence by Corollary \ref{Saturated1}, the ideal
$(I_{Y_1})_{\leq d}$ is a saturated ideal defining a scheme
consisting of an unmixed curve $C_1$ of degree $s$ plus some zero-dimensional
scheme. But $C$ and $C_1$ have infinitely many hyperplane sections
which agree. This means $C=C_1$ since they are unmixed. Now we
claim that the ideal $\langle(I_{Y_1})_{\leq d}\rangle$ is precisely
the ideal of $C$. We know that $\langle(I_{Y_1})_{\leq d}\rangle$
coincides with $I_C$ in degree $\leq d$. Hence if we prove that $C$
is $d$-regular then we will have proved our claim. But this follows from
Proposition \ref{firstCH regular} since
$d>r_2(R/I_C)=r_2(R/I_{Y_1})$ and (4.7) holds. This proves (c), the
first half and second part of (e). Now we will prove (f). Suppose
that $h^1(\mathcal I_{C_\red}(d-1))=0$. Note that we know
$I_C\subset I_{C_{\red}}\subset I_{Y_1}$ since $Y_1$ is reduced,
so we have that $(I_C)_t=(I_{C_{\red}})_t$ for all $t\leq d+1$.
Hence $d>r_2(R/I_{C_{\red}})=r_2(R/I_C)$. This implies $I_{C_{\red}}$
is $d$-regular (Proposition \ref{firstCH regular}) and we get
$C=C_{\red}$. This means the scheme $V$ defined by the saturated
ideal $\langle(I_Y)_{\leq d}\rangle$ is reduced (the proof is
exactly same as that of \cite{BGM}).
\end{proof}

\begin{Remk} \label{compare BGM}
Theorem \ref{Generalize} generalizes the analogous result (Corollary 3.7)  in
\linebreak
\cite{BGM}: we remove the condition $d\geq s$.  The next result is the main
consequence for us, and introduces the key study of seeing what happens when we
replace ``$d \geq s$" with
``$d > r_2(R/I_Y)$.''  It generalizes \cite{BGM} Theorem 3.6. The
bulk of this study
will be in the next section when we consider uniform position, UPP.  For now,
though, we note that the only conclusion of \cite{BGM} Theorem 3.6 that is not
included in our Theorem \ref{Generalize} is that we do not claim that $V$ is
reduced, except in part (f), with an extra hypothesis.  In the next
section we will
see that this is not necessarily true without that hypothesis.
\end{Remk}

\begin{Coro}\label{Coro_main}
Let $Y\subset \mathbb P^{n-1}$, $n\geq 3$, be a closed subscheme of
dimension $\leq 1$ in $\mathbb P^{n-1}$. Assume that
\[\Delta H(Y,d)= \Delta H(Y,d+1)=s\neq 0\]
for some $d> r_2(R/I_Y)$. Then (a)-(f) of Theorem {\rm
\ref{Generalize}} continue to hold.
\end{Coro}

\begin{proof}
This follows from Corollary \ref{Saturated} and Theorem
\ref{Generalize}.
\end{proof}


\section{Application to points with UPP} \label{UPP sect}

In this section we add the hypothesis that our finite set of points has the
Uniform Position Property (UPP).  The following was shown in Theorem 4.7 of
\cite{BGM}.

\begin{Thm} \label{BGM result quote}
Let $Z \subset {\mathbb P}^{r+1}$ be a reduced finite set of points with UPP.
Assume that $\Delta H(R/I_Z,d) = \Delta H(R/I_Z,d+1) = s$, $d \geq
s$.  Then there
exists a reduced, irreducible curve $C$ of degree $s$ such that

\begin{itemize}
\item[(a)] $Z \subset C$;

\item[(b)] $I_C = \langle (I_Z)_{\leq d} \rangle$;

\item[(c)] $\Delta H(R/I_Z,t) = \Delta H(R/I_C,t)$ for all $t \leq d+1$.
\end{itemize}
\end{Thm}

In this section we assume that $Z$ has UPP, and we see what
consequences this has
under our hypothesis $d > r_2(R/I_Z)$.  We also give several examples of
``expected'' behavior that does {\em not} occur.  In doing so, we give some new
insight into what can be the Hilbert function of a set of points
$Y$ with UPP.  We stress again  that this hypothesis differs from the
analogous one
in
\cite{BGM} only in that there it was assumed that $d\geq s$.

We begin this section by seeing what is still true for
points with UPP. The rest of this section focuses on what is no longer true
(even if counterexamples are not so easy to come by).

\begin{Thm}\label{UPP}
Let $Z \subset \mathbb P^{n-1}$ be a reduced finite set of
points with UPP and let $(1,h_1,\ldots,h_t)$ be the $h$-vector of
$Z$. Assume that
\[h_d=h_{d+1}=s\neq 0 \]
for some $d>r_2(R/I_{Z})$. Then there exists an unmixed  curve $C$ of
degree $s$
such that
\begin{enumerate}
\item[(a)] $Z \subset C$

\item[(b)] $I_C=\langle(I_{Z})_{\leq d}\rangle$

\item[(c)] $\Delta H(C,t)=h_t$ for all $t\leq d+1$.
\end{enumerate}
\end{Thm}

\begin{proof}
By Theorem \ref{Generalize} and Corollary \ref{Coro_main}, we know that there
exists an unmixed curve $C$ of degree $s$ containing a subset
${Z}_1\subset Z$, such
that $I_{{Z}_1}$ agrees with $I_C$ up to degree $d+1$. Since
$Z$ has UPP, $(I_{Z})_t=(I_{Z_1})_t$ for $t\leq
d+1$. Hence $I_Z$, $I_{Z_1}$ and $I_C$ all agree up to degree $d+1$.
On the other hand,  Corollary \ref{Coro_main} (e) says that
$\Delta H(R/I_{Z_1},t)=\Delta H(R/I_{Z},t)$ for $t\geq d$, so we get that
$Z=Z_1$. This proves (a), (b) and (c).
\end{proof}

Apart from the hypotheses on $d$, the only difference between Theorem \ref{BGM
result quote} and Theorem \ref{UPP} is that in the former (with the stronger
hypothesis $d \geq s$), it was possible to conclude that $C$ was reduced and
irreducible.

Another property that one might expect from a set of points with UPP
is the following.  Suppose that $\Delta H(R/I_Z,d) = \Delta H(R/I_Z,d+1) = s
>  \Delta H(R/I_Z,d+2)$ for some   $d > r_2(R/I_Z)$.  Then one might expect
that
the first difference $\Delta H(R/I_Z)$ is strictly decreasing from that point,
i.e.\
$\Delta H(R/I_Z,t) > \Delta H(R/I_Z,t+1)$ for all $t \geq d+1$, until
it becomes
zero.  Corollary
\ref{str dec UPP} shows that in fact this is true if $d \geq s$, the
assumption from
\cite{BGM} (although it was not proved there), and Proposition \ref{strictly
decreasing} gives a more general condition that guarantees this property.

\begin{Prop} \label{strictly decreasing}
Let $Z$ be a zero-dimensional scheme with
\[
s = \Delta H(R/I_Z,d) = \Delta H(R/I_Z,d+1) > \Delta H(R/I_Z,d+2)
\]
for $d > r_2(R/I_Z)$, and assume that the unmixed curve $C$ of degree
$s$ guaranteed
by Corollary \ref{Coro_main} is reduced and irreducible and contains
all of $Z$.
Then
\begin{equation} \label{str gr}
\Delta H(R/I_Z,t) > \Delta H(R/I_Z,t+1) \ \hbox{ for $t \geq d+1$  as long as
$\Delta H(R/I_Z,t) > 0$.}
\end{equation}

In particular, suppose that $Z$ is a reduced set of points with UPP.
Then if $C$
is reduced and irredcible, we get (\ref{str gr}).
\end{Prop}

\begin{proof}
If it does not hold then say $\Delta H(R/I_Z,e) = \Delta H(R/I_Z,e+1)
= s' > 0$ for
some $e \geq d+2$.  We already know that $s' < s$, thanks to Proposition
\ref{Decreasing} and the hypothesis $\Delta H(R/I_Z,d+1) > \Delta
H(R/I_Z,d+2)$.  We
obtain a curve $C'$ of degree $s'$, thanks to Theorem \ref{Saturated1},  and
clearly $C'$ is contained in $C$.  But since $C$ is reduced and
irreducible, this is
impossible.
\end{proof}

\begin{Coro} \label{str dec UPP}
If $Z$ is a reduced set of points with UPP, and
\[
s = \Delta H(R/I_Z,d) = \Delta H(R/I_Z,d+1) > \Delta H(R/I_Z,d+2)
\]
for $d \geq s$ then (\ref{str gr}) holds.
\end{Coro}

\begin{proof}
This follows from \cite{BGM} Theorem 4.7, and Proposition
\ref{strictly decreasing}
above.
\end{proof}

\begin{Remk}
Notice again that Theorem \ref{UPP} does not claim that $C$ is reduced
or irreducible, even though these conclusions do appear in \cite{BGM} in the
analogous results, under the stronger assumption that $d \geq s$.  In
the examples
that follow, we will see that these conclusions may in fact be false
under only the
assumption $d > r_2(R/I_Y)$!  This is surprising since prior to this
our results
matched those of \cite{BGM} very closely, and since UPP is strongly tied to
irreducibility in many ways (see for instance the discussion in the
introduction).
\end{Remk}

\begin{Ex} \label{CD example}
We first recall a beautiful example from \cite{CD}, which will set
the stage for
our first main example.  Consider the family of complete intersection ideals
$I_{m,n} := (x^mt - y^mz, z^{n+2}-xt^{n+1})
\subset k[x,y,z,t]$.  Then for
$m,n \geq 1$, $\hbox{reg}(I_{m,n}) = m+n+2$ while $\hbox{reg}(\sqrt{I_{m,n}}) =
mn+2$.  This settles (in an unexpected way) a long-standing and very
interesting
question (cf.\
\cite{ravi}) of whether it is possible that the regularity increases if one
replaces an ideal by its radical.

Now, we take $m = n = 4$.  The complete intersection $I_{4,4}$ is then of type
$(5,6)$, has degree 30, and has a Hilbert function whose first difference is

\begin{center}
\begin{tabular}{l|cccccccccccccccccccccccccccccc}
deg & 0 & 1 & 2 & 3 & 4 & 5 & 6 & 7 & 8 & 9 & 10 & 11 & 12 \\
$\Delta H$ & 1 & 3 & 6 & 10 & 15 & 20 & 24 & 27 & 29 & 30 & 30 & \dots
\end{tabular}

\end{center}

On the other hand, $\sqrt{I_{4,4}}$ can be computed on a computer
program, e.g.\
{\tt macaulay} \cite{macaulay}.  It has degree 26 (hence $I_{4,4}$ is not
reduced), and its Hilbert function has first difference

\medskip

{\small
\begin{center}
\begin{tabular}{l|cccccccccccccccccccccccccccccc}
deg & 0 & $\dots$ & 6 & 7 & 8 & 9 & 10 & 11 & 12 & 13 & 14 & 15 & 16 &
17 & 18 & 19 & 20
\\ \hline
$\Delta H$ & 1 & $\dots$  & 24 & 27 & 29 & 29 & 29 & 29 & 28 & 28 &
28 & 27 & 27 & 27 & 26 & 26 &
\dots
\end{tabular}

\end{center}
}
\noindent (where the early degrees agree with $I_{4,4}$).

\medskip

We focus on  the ideal $ \sqrt{I_{4,4}}$.  Of course it defines a reduced
subscheme of ${\mathbb P}^3$.  One can verify on the computer (but see below as
well) that
$r_2(R/\sqrt{I_{4,4}}) = 8$, as one would expect.  Setting $d = 10$,
though, one
sees that $d > r_2(\sqrt{I_{4,4}})$ and $\Delta H(R/\sqrt{I_{4,4}}),d) = \Delta
H(R/\sqrt{I_{4,4}}) = 29$, so $\langle (\sqrt{I_{4,4}})_{\leq 10}
\rangle$ is the
saturated ideal of a curve of degree 29 (Theorem \ref{Generalize}
above).  One can
again check that this curve is not reduced, by checking that its radical is
precisely $\sqrt{I_{4,4}}$ (which has degree 26), so one also sees
that it is not
irreducible (by degree considerations).
\end{Ex}

Example \ref{CD example} does not involve a finite set of points.
However, we will
now give a refinement of this (Example \ref{CD gen}).  This example,
together with
Example \ref{nonreduced} and Example \ref{not irred}, show some interesting and
unexpected possible behavior of points with UPP.  Contrast them with Theorem
\ref{BGM result quote} above.

\begin{Ex} \label{CD gen}
We again work in the ring $R = k[x,y,z,t]$ and consider the ideal $I_\lambda =
(z,t)$.  This defines a line, $\lambda$, in ${\mathbb P}^3$.  Let $F \in
I_\lambda$ be a homogeneous polynomial of degree 5 defining a smooth surface
containing
$\lambda$ (in particular, it is smooth at all points of
$\lambda$).  By slight abuse of notation, we will use $F$ also for the surface.
Consider the ideal $I' = I_\lambda^5 + (F)$.  It can be checked from the exact
sequence
\[
0 \rightarrow I_\lambda^4(-5) \rightarrow I_\lambda^5 \oplus R(-5) \rightarrow
I_\lambda^5 + (F)
\rightarrow 0
\]
(or by direct computation on the computer)
that $I'$ is the saturated ideal of a curve (not  arithmetically
Cohen-Macaulay) of degree 5, corresponding to the divisor $D :=
5\lambda$ on the
smooth surface $F$.  Consider the linear system $|6H-D|$ on $F$.  We
first claim
that this linear system has no base locus.  Indeed, consider first the curve
$\lambda$ and the linear system $|H-\lambda|$.  By considering the pencil of
planes through $\lambda$ and the residual cut out on $F$, we see that
this linear
system has no base locus (otherwise a point $P$ of the base locus would be a
singular point of $F$).  Hence
$|5H-D|$ has no base locus, since a union of five planes containing $\lambda$
also contains $D$, so all the more $|6H-D|$ has no base locus.

By Bertini's theorem, then, the general element $C$ of $|6H-D|$ is smooth.
Notice that the ideal $I'$ contains sextics in addition to the elements of
$(I_\lambda^5)_6$, so the general element of $|6H-D|$ is also irreducible.

What is this curve $C$?  It is simply the residual to $D$
(viewed now as a curve in
${\mathbb P}^3$) in the complete intersection of $F$ and a general
element $G \in
I'$ of degree 6.  So this residual has degree $30-5 = 25$.  Its
Hilbert function
turns out to have first difference

{\scriptsize
\begin{center}
\begin{tabular}{l|cccccccccccccccccccccccccccccc}
deg & 0 & $\dots$ & 7 & 8 & 9 & 10 & 11 & 12 & 13 & 14 & 15 & 16 &
17 & 18 & 19 & 20 & 21 & 22
\\ \hline
$\Delta H$ & 1 & $\dots$ & 27 & 29 & 29 & 29 & 29 & 28 & 28 &
28 & 27 & 27 & 27 & 26 & 26 & 26 & 25 & \dots
\end{tabular}
\end{center}
}
\medskip
\noindent (the entries of low degree are the same as in Example
\ref{CD example};
all the values after degree 21 are 25).
Note that this Hilbert function agrees with that of the example of Chardin and
D'Cruz (Example \ref{CD example}) up to and including degree 20.  What we know
that is new (and not true in their example) is that $C$ is smooth.

Now let $Z$ consist of sufficiently many generally chosen points on
$C$.  $Z$ has
UPP since $C$ is smooth and irreducible, and by choosing enough points we can
assume that the Hilbert function of $Z$ agrees with that of $C$ past degree 21.
Then $r_2(R/I_Z) = 8$ and $\Delta H(R/I_Z,10) = \Delta H(R/I_Z,11) = 29$.  The
ideal generated by the components of degrees $\leq 10$ then defines a curve of
degree 29.  Since this ideal is contained in $I_Z$, this curve
contains $C$ (and
hence also $Z$).  But in fact this curve of degree 29 consists of the
residual in
the complete intersection of degree 30 to the line $\lambda$, hence it consists
of the union of $C$ and a curve of degree 4 supported on $\lambda$.
So this curve
defined by the degree 10 component of $I_Z$ is neither reduced nor
irreducible, even
though $Z$ has UPP.

To summarize, this example  shows that a set of points $Z$ with UPP, for which
\[
\Delta H(R/I_Z,d) =
\Delta H(R/I_Z,d+1) = s > \Delta H(R/I_Z,d+2) \hbox{ for $d > r_2(R/I_Z)$},
\]
  can have the
component of degree
$d$ define a non-reduced, non-irreducible curve (but it has to have
degree $s$),
and the Hilbert function can fail to be strictly decreasing past degree $d+2$.
\end{Ex}

One might wonder if Example \ref{CD gen} was somehow an ``accident,"
in the sense
that while the component of $I_Z$ in question was neither reduced nor
irreducible,
still all the points of $Z$ lay one one irreducible and reduced component.
We now make a minor adjustment to show that even this is not necessarily true.

\begin{Ex} \label{nonreduced}
In Example \ref{CD gen}, instead of choosing ``sufficiently many"
points on $C$,
instead choose $Z$ to consist of 192 general points on $C$ and one
general point of
$\lambda$.  The first difference of the Hilbert function of $Z$ is

\begin{center}
\begin{tabular}{l|cccccccccccccccccccccccccccccc}
deg & 0 & 1 & 2 & 3 & 4 & 5 & 6 & 7 & 8 & 9 & 10 & 11 &   \\
\hline $\Delta H$ & 1 & 3 & 6 & 10 & 15 & 20 & 24 & 27 & 29 & 29 &
29 & 0
\end{tabular}

\end{center}

Note that this is exactly the same as what we would have if we had taken 193
general points of $C$.  Again, $r_2(R/I_Z) = 8$.  The base locus of
$(I_C)_{9}$ and
$(I_C)_{10}$ is exactly the non-reduced and reducible curve of degree
29 mentioned
in Example~\ref{CD gen}, so the Hilbert function of $Z$ is the
truncation of the
Hilbert function given above.  And by the general choice of the
points, this will
continue to be true regardless of which subsets we take.  Hence $Z$ has UPP and
satisfies $\Delta H(R/I_Z,d) = \Delta H(R/I_Z,d+1) = s$ for some $d >
r_2(R/I_Z)$,
but not all of the points of $Z$ lie on one reduced and irreducible
component of the
curve of degree $s$ obtained by our result (since one point lies on
the non-reduced
component).
\end{Ex}

Again, one may wonder if it could happen that at least ``almost all'' of the
points must lie on one reduced and irreducible component.  We now
give an example
where the entire curve arising from Theorem
\ref{UPP} is reduced, but consists of two irreducible components, and
each of these
contains half of the points; but yet the points still have UPP.

\begin{Ex} \label{not irred}
As remarked above, it is shown in \cite{BGM} that if $\Gamma$ has  UPP and if
$\Delta H(R/I_\Gamma,d) = \Delta H(R/I_\Gamma,d+1) = s$ for $d \geq s$ then
$\Gamma$ lies on a reduced and irreducible curve of degree $s$.  We now
show that with only the assumption $d > r_2(R/I_\Gamma)$, this curve
need not be
irreducible (although in Theorem \ref{UPP} above we did prove the
existence of this
curve).

Let $Q$ be a smooth quadric surface in $\mathbb P^3$, and by abuse of
notation we
use the same letter $Q$ to denote the quadratic form defining $Q$.
Let $C_1$ be a
{\em general} curve on $Q$ of type $(1,15)$, and let $C_2$ be a {\em general}
curve on $Q$ of type $(15,1)$.  Hence both $C_1$ and $C_2$ are smooth rational
curves of degree 16, and $C := C_1 \cup C_2$ is the complete
intersection of $Q$
and a form of degree 16.  Note that $C$ is arithmetically Cohen-Macaulay, but
$C_1$ and $C_2$ are not.

It is not difficult to compute the Hilbert functions of these curves.
We record
their first differences (of course there is no difference between behavior of
$C_1$ and behavior of $C_2$; this will be important in our argument):

\bigskip

{\small

\begin{tabular}{l|ccccccccccccccccccccccccccccccccccccccccccc}
degree & 0 & 1 & $\dots$ & 7 & 8 & 9 & 10 & 11 & 12 & 13 & 14 & 15 & 16
& 17 & 18 \\ \hline
$\Delta H(R/I_C,-)$ & 1 & 3 & $\dots$ & 15 & 17 & 19 & 21 & 23 & 25 &
27 & 29 & 31 & 32 & 32 & 32 \\
$\Delta H(R/I_{C_i},-)$ & 1 & 3 & $\dots$ & 15 & 17 & 19 & 21 & 23 &
25 & 27 & 29 & 16 & 16 & 16 & 16 \\
\end{tabular} }

\medskip

\noindent We now observe:
\begin{enumerate}

\item These first differences (hence
the ideals themselves) agree through degree 14, and in fact the only generator
before degree 15 is $Q$.

\item \label{hf values} By adding these values, we see that $H(R/I_C,18) = 352$
and
$H(R/I_{C_i},15) = 241$.

\item Since $C$ and $C_i$ are curves, these values represent the Hilbert
functions of
$I_C + (L)$ and $I_{C_i} +(L)$ for a general linear form $L$.

\item \label{value of r2} $r_2(R/I_C) = 16$ since $C$ is an arithmetically
Cohen-Macaulay curve.
\end{enumerate}

\bigskip

Let $\Gamma_1$ (respectively $\Gamma_2$) be a general set of 176
points on $C_1$
(respectively $C_2$).  So $\Gamma := \Gamma_1 \cup \Gamma_2$ is a set of 352
points whose Hilbert function agrees with that of $C$ through degree 18.  In
particular, we have $\Delta H(R/I_\Gamma,17) = \Delta
H(R/I_\Gamma,18) = 32 = \deg
C$.  Furthermore, $r_2(R/I_\Gamma) = 16$ by (\ref{value of r2}) above.  Hence
Theorem \ref{UPP} applies, and we indeed have  that the component of
$I_\Gamma$ in
degree 17 defines $C$.  However, $C$ is not irreducible.  Our goal is to show
that $\Gamma$ has the Uniform Position Property, and thus there is no chance of
showing that all the points must lie on a unique irreducible
component in Theorem
\ref{UPP}, under our hypothesis (as was done in
\cite{BGM}).

To show UPP, it is enough to show that the union of any choice of
$t_1$ points of
$\Gamma_1$ (i.e.\ $t_1$ general points of $C_1$ for $t_1 \leq 176$) and $t_2$
points of $\Gamma_2$ (i.e.\ $t_2$ general points of $C_2$ for $t_2
\leq 176$) has
the truncated Hilbert function.  For example, if $t_1 = 150$ and $t_2
= 160$ then
we have to show that
$\Delta H(R/I_\Gamma)$ has values
\[
1 \ \ 3 \ \ 5 \ \ 7 \ \ 9 \ \ 11 \ \ 13 \ \ 15 \ \ 17 \ \ 19 \ \ 21 \
\ 23 \ \ 25
\ \ 27 \ \ 29 \ \ 31 \ \ 32 \ \ 22 \  \ 0.
\]
Notice that we know that some subset has this Hilbert function, by \cite{GMR}.
We have to show that {\em all} subsets have this Hilbert function.

Let $Y = Y_1 \cup Y_2$ be our choice of $t_1$ points of $\Gamma_1$ and
$t_2$ points of $\Gamma_2$.  Without loss of generality, we may assume that
$t_1 \geq t_2$.  Hence because of the indistinguishability of $C_1$
and $C_2$ and
the generality of the points, for any degree $d$ it may be that every
element of
$(I_Y)_d$ vanishes on all of $C = C_1 \cup C_2$, or it may be that
every element
of
$(I_Y)_d$ vanishes identically on
$C_1$ but not $C_2$ (since $t_1 \geq t_2$), but (*) {\em it cannot happen that
every element of
$(I_Y)_d$ vanishes identically on $C_2$ but not on all of $C_1$.}

Now, we have to determine the Hilbert function of $R/I_Y$.  {\em The Hilbert
function of the union does not depend on what order we consider these points. }
First consider $Y_1$.  Since $Y_1$ consists of $t_1$ general points
of $C_1$,  the
Hilbert function of $R/I_{Y_1}$ is the truncation of that of
$R/I_{C_1}$.  Thanks
to item (\ref{hf values}.) above, and since $t_1 \leq 176 < 241$, the Hilbert
function of $R/I_{Y_1}$ is the truncation of $R/I_C$ as well.  We now add the
points of $Y_2$ one by one.

Suppose that at some point, adding a point $P$ of $Y_2$ does not result in a
Hilbert function that is a truncation of $R/I_C$.  Suppose that previous to $P$
we have gone through a subset $Y_2' \subset Y_2$ and gotten a truncated Hilbert
function each time, so $P$ is the first time that this fails.  This means that
there is some degree $d$ such that
\[
\Delta H(R/I_C,d) > \Delta H(R/I_{Y_1 \cup Y_2'},d) = \Delta
H(R/I_{Y_1 \cup Y_2'
\cup P},d).
\]
In other words, every form of degree $d$ containing $Y_1 \cup Y_2'$
also contains
$P$, but there is some form $F$ of degree $d$ containing $Y_1 \cup
Y_2'$ but not
all of $C$.  Now, $P$ is a general point of $C_2$ independent of the choices of
points previous to it, and every element of $(I_{Y_1 \cup Y_2'})_d$ vanishes at
$P$.  Therefore every element of $(I_{Y_1 \cup Y_2'})_d$ vanishes on an open
subset of $C_2$, and hence (by irreducibility) on all of $C_2$.  Since there is
an $F \in (I_{Y_1 \cup Y_2'})_d$ not vanishing on all of $C$, $F$
must not vanish
on all of $C_1$.  But our conclusion (*) above used only the fact that $t_2
\leq t_1$, so it remains true replacing $Y_2$ by $Y_2' \cup P$.  Contradiction.
Therefore $\Gamma$ has UPP.
\end{Ex}


\section{The connection to the Weak Lefschetz Property} \label{add'l comments}

In the previous sections, we gave the results on the Hilbert
function of the scheme $Y$ of $\dim(Y)\leq 1$ in degree
$d>r_2(R/I_Y)$ in terms of its first difference.  The range in which we obtained
our results was $d > r_2(R/I_Y)$.  We also gave examples to show that we cannot
expect these results to extend to the range $d \leq r_2(R/I_Y)$ (Remark \ref{not
curve}).  However, in this section we show that in fact there are useful
analogous results that we can obtain in this smaller range, but changing the
focus somewhat (Proposition
\ref{Decreasing-1} through Theorem
\ref{r2 case}).  As an application, we give new results  in terms of second
difference of the Hilbert function if $Y$ is a set of points with WLP in
$\mathbb P^{n-1}$,
$n>3$. In particular, we apply these result to a set of points with UPP and WLP
in
$\mathbb P^3$ (Theorem
\ref{UPP_in_P3} and Corollary \ref{UPP_WLP}).  It is a natural question whether
every set of points with UPP has WLP, since both are open properties.  (It is
easy to find examples having WLP but not UPP.)  We give an example (Example
\ref{UPP vs WLP}) to show that this is not the case.

First we prove analogies to Proposition \ref{Decreasing} and
Theorem \ref{Saturated1} for $d\leq r_2(R/I_Y)$.

\begin{Prop}\label{Decreasing-1}
Let $Y$ be a scheme of dimension $\leq 1$ in $\mathbb P^{n-1}$,
$n>3$. Let $K = {(I_Y + (L_1,L_2))}/{(L_1,L_2)}$; ;this is a homogeneous
ideal in $S=R/(L_1,L_2)$. Then
\begin{equation}\label{eq6}
H(S/K,d)\geq H(S/K,d+1)
\end{equation}
for $d\geq r_3(R/I_Y)$. Moreover, for $r_2(R/I_Y)>d>r_3(R/I_{Y})$,
we have equality in {\rm({\ref{eq6}})} if and only if $\Gin (K)$
has no minimal generators in degree $d+1$.
\end{Prop}
\begin{proof}
The proof can be given in the same manner as that of Proposition
\ref{Decreasing} (though not trivially).
\end{proof}

\begin{Prop}\label{Saturated-2}
Under the same situation as Proposition \ref{Decreasing-1}, if we
have equality in {\rm({\ref{eq6}})}, say $s$, for
$r_2(R/I_Y)>d>r_3(R/I_{Y})$ then $\langle (I_Y)_{\leq
d}\rangle^{\sat}$ is a homogeneous ideal defining a
two dimensional subscheme of degree $s$ in $\mathbb P^{n-1}$.
\end{Prop}
\begin{proof}
If we consider the homogeneous ideal of $S=R/(L_1,L_2)$
\[M = \frac{(\langle (I_Y)_{\leq d}\rangle + (L_1,L_2))}{(L_1,L_2)}\]
then $M$ is $d$-regular by Theorem \ref{CP} and Proposition
\ref{Decreasing-1}. Consider the generic initial ideal of $M$.
Since $M$ is $d$-regular for $r_3(R/I_Y)<d<r_2(R/I_Y)$, $M$ could
not have any power of $x_{n-2}$.  This means that $D(M)=n-3$,
since $d>r_3(R/I_Y)$ implies that the degree $d$ component of $Gin(I_Y)$ contains
the monomial $x_{n-3}^d$. Note that
$\Gin(M)=(\Gin(\langle (I_Y)_{\leq d}\rangle)_{x_n\rightarrow
0})_{x_{n-1}\rightarrow 0}$ does not have any minimal generators of
degree larger than $d$. Hence every minimal generator of $\langle
(I_Y)_{\leq d}\rangle$ in degree $>d$ must involve the variables
$x_{n-1}$ or $x_n$. This implies $D(\langle (I_Y)_{\leq
d}\rangle)=n-3$. Hence $\langle (I_Y)_{\leq d}\rangle^{\sat}$
defines a two dimensional subscheme of $\mathbb P^{n-1}$.

It follows from $\Delta^2H(R/\langle (I_Y)_{\leq
d}\rangle,t)=H(S/M,t)=s$ for sufficiently large $t$ that the
degree of $(R/\langle (I_Y)_{\leq d}\rangle^{\sat})$ is exactly
$s$.
\end{proof}

\begin{Coro}\label{Common factor}
Let $Y$ be a scheme of dimension $\leq 1$ in $\mathbb P^{3}$. If
we have equality in {\rm({\ref{eq6}})}, say $s$, for
$r_2(R/I_Y)>d>r_3(R/I_{Y})$, then $(I_Y)_d$ has a common factor
$F$ of degree $ s$.
\end{Coro}
\begin{proof}
Let $X$ be the unmixed part of $\langle (I_Y)_{\leq d}\rangle$ of
$\dim(X)=2$. Then $I_{X}=(F)$ for a  homogeneous
polynomial $F$ since $X$ has codimension $1$ in $\mathbb P^3$.
Hence $(I_Y)_d\subset (I_X)_d\subset(F)_d$ and thus $(I_Y)_d$ has
a common factor $F$ of degree $ s$.
\end{proof}

We need an additional condition to be able to deduce the saturatedness and
regularity of $\langle (I_Y)_{\leq d}\rangle$.

\begin{Thm}\label{r2 case}
Let $Y$ be a scheme of dimension $\leq 1$ in $\mathbb P^{n-1}$,
$n>3$. Let $K = {(I_Y + (L_1,L_2))}/{(L_1,L_2)} \subset S=R/(L_1,L_2)$. Suppose
that $H(S/K,d) =$ \linebreak $ H(S/K,d+1)$ for $r_2(R/I_Y)>d>r_3(R/I_{Y})$.  
Then $\langle (I_Y)_{\leq d}\rangle$ is a homogeneous ideal (not necessarily
saturated) defining a two dimensional scheme of degree $s$ in $\mathbb
P^{n-1}$.
Furthermore, the
following are  equivalent.
\begin{enumerate}
\item[(a)] $\langle (I_Y)_{\leq d}\rangle$ is  $d$-regular.
\item[(b)] $0 \rightarrow
(R/I_Y+(L_1))_d \stackrel{\times L_2}{\longrightarrow}
(R/I_Y+(L_1))_{d+1}$ for a general linear form $L_2$.
\end{enumerate}
If these equivalent conditions are satisfied then the ideal is
saturated as well.
\end{Thm}
\begin{proof}
The first assertion is just a repetition of Proposition \ref{Saturated-2}, just
for completeness. We now prove the equivalence. Let $\bar I=\langle (I_Y)_{\leq
d}\rangle$. For a general linear form $L_2$ if the map
\[
(R/I_Y+(L_1))_d \stackrel{\times L_2}{\longrightarrow}
(R/I_Y+(L_1))_{d+1}
\] 
is not injective, then we know that there is
a minimal generator of $\Gin(\bar I)$ with degree $d+1$ from the
proof of Proposition \ref{Decreasing}. But this means $\langle
(I_Y)_{\leq d}\rangle$ is not $d$-regular (Theorem \ref{BS_reg}).

Conversely, suppose that the map
\[
(R/I_Y+(L_1))_d \stackrel{\times L_2}{\longrightarrow}
(R/I_Y+(L_1))_{d+1}
\] 
is injective for a general linear form
$L_2$. It is enough to show that $\Gin(\bar  I)$ has no minimal
generator of degree $d+1$ in order to prove $\langle (I_Y)_{\leq
d}\rangle$ is $d$-regular. Note that $(\Gin(\bar
I)_{x_{n}\rightarrow 0})_{x_{n-1}\rightarrow 0}$ is $d$-regular.
Hence if $\Gin(\bar  I)_{x_{n}\rightarrow 0}$ has a minimal
generator of degree $d+1$ then it must involve the variable
$x_{n-1}$. But this is impossible because the map $(R/I_Y+(L_1))_d
{\longrightarrow} (R/I_Y+(L_1))_{d+1}$ is injective. Hence
$\Gin(\bar  I)_{x_{n}\rightarrow 0}$ is $d$-regular (Theorem
\ref{CP}). Suppose that there is a minimal generator $x^L$ of
$\Gin(\bar  I)$ with degree $d+1$. Then it must involve the
variable $x_n$. Notice that $\Gin(\bar
I)_{d+1}\subset\Gin(I_Y)_{d+1}$ and $x^L$ cannot be a minimal
generator of $\Gin(I_Y)$ since $I_Y$ is a saturated ideal. Hence
we can choose a monomial $x^T \in \Gin(I_Y)_{d}$ such that
$x_ix^T=x^L$ for some variable $x_i$. But it contradicts the fact
that $x^L$ is a minimal generator of $\Gin(\bar  I)$ because
$x^T\in \Gin(\bar  I)_{d}=\Gin(I_Y)_{d}$. Hence $\bar  I$ is
$d$-regular (Theorem \ref{CP}). 

The proof that then $\bar  I$ is
saturated is exactly same as that of Theorem \ref{Saturated1}
\end{proof}

\begin{Remk}
The proof of Theorem \ref{r2 case} makes it clear that if 
\[
(R/I_Y+(L_1))_d \stackrel{\times L_2}{\longrightarrow}
(R/I_Y+(L_1))_{d+1}
\]
 is not injective then $(I_Y)_{\leq d}$ is not $d$-regular.  However, the
proof leaves it unclear whether $(I_Y)_{\leq d}$ is saturated in this case.  In
the following example, we show that this ideal may or may not be saturated.  We
used {\tt macaulay} for this calculation, but it could also be computed by
hand.
\end{Remk}

\begin{Ex} \label{not having WLP}
Let $J = (F_1,F_2,F_3) \subset R$ be a regular sequence of
type $(2,2,2)$ defining a zero-dimensional complete intersection $Z_1$ of degree
8 in $\mathbb P^3$, and let $Q$ be a general quadratic polynomial.  (It is enough
that $Q$ not vanish on any of the 8 points of $Z$.  We will also use $Q$ to
denote the corresponding quadric surface.)  Let $I = (QF_1,QF_2,QF_3)$.  This
is the saturated ideal of the union of $Z_1$ and $Q$.  Let
$Z_2$ be a set of 81 general points on $Q$, and let $Y = Z_1 \cup
Z_2$.  As above, let $L_1$ and $L_2$ be general linear forms.  The Hilbert
functions of the relevant ideals are as follows.

\bigskip

\begin{tabular}{l|ccccccccccccccccccccccccccccccccccccccccccc}
degree & 0 & 1 & 2 & 3 & 4 & 5 & 6 & 7 & 8 & 9 & $\dots$ \\ \hline
$H(R/I_Y,-)$ & 1 & 4 & 10 & 20 & 32 & 44 & 57 & 72 & 89 & 89 & $\dots$  \\
$H(R/(I_Y+(L_1)),-)$ & 1 & 3 & 6 & 10 & 12 & 12 & 13 & 15 & 17 & 0 &  $\dots$
\\
$H(R/(I_Y+(L_1,L_2)),-)$ & 1 & 2 & 3 & 4 & 2 & 2 & 2 & 2 & 2 & 0 & $\dots$
\end{tabular} 

\bigskip

\noindent From this we see that $r_2(R/I_Y) = 8$, and it is not hard to see
that $r_3(R/I_Y) = 3$. 

Let $d = 4$.   We have the needed equality in Corollary
\ref{Common factor}, so $(I_Y)_4$ has a common factor of degree 2, namely $Q$. 
In fact,
$(I_Y)_{\leq 4}$ is exactly $I$, and it is not hard to check that its
regularity is 6 (since the regularity of $(F_1,F_2,F_3)$ is 4).  Since 
\linebreak $H(R/(I_Y+(L_1,L_2)),5)= 2$, we see that we do not have the
injectivity of the desired map.  Hence we confirm in this example what is
proved in Theorem
\ref{r2 case}, that we do not have regularity in degree $d =4$.  However, as
mentioned above, $(I_Y)_{\leq 4} = I$ is saturated.  Note that the
non-injectivity implies that
$Y$ does not have the Weak Lefschetz Property.

Now we modify this example a little bit.  Let $Z_1$ now consist of 16 general
points in $\mathbb P^3$ and let $Q$ again be a quadric surface not containing
any of the points of $Z_1$.  Note that  $I_{Z_1}$ has four cubic generators
and three quartic generators, but the component in degree 3 does define $Z_1$
scheme-theoretically.  Let $Z_2$ again consist of 81 general points of $Q$,
and let $Y = Z_1 \cup Z_2$.  We have the following Hilbert functions.

\bigskip

\begin{tabular}{l|ccccccccccccccccccccccccccccccccccccccccccc}
degree & 0 & 1 & 2 & 3 & 4 & 5 & 6 & 7 & 8 & 9 & $\dots$ \\ \hline
$H(R/I_Y,-)$ & 1 & 4 & 10 & 20 & 35 & 52 & 65 & 80 & 97 & 97 & $\dots$  \\
$H(R/(I_Y+(L_1)),-)$ & 1 & 3 & 6 & 10 & 15 & 17 & 13 & 15 & 17 & 0 &  $\dots$
\\
$H(R/(I_Y+(L_1,L_2)),-)$ & 1 & 2 & 3 & 4 & 5 & 2 & 2 & 2 & 2 & 0 & $\dots$
\end{tabular} 

\bigskip

\noindent This time $r_2(R/I_Y)$ is again 8, but $r_3(R/I_Y) = 4$.  Let $d =
5$.  We again have the needed equality, so $(I_Y)_5$ has a common factor of
degree 2, namely $Q$.  

The map
\[
(R/I_Y+(L_1))_5 \stackrel{\times L_2}{\longrightarrow}
(R/I_Y+(L_1))_{6}
\]
obviously has no chance to be injective.  Hence by Theorem \ref{r2 case},
$(I_Y)_{\leq 5}$ is not 5-regular.  What about saturation?  Notice that
$(I_Y)_{\leq 5}$ is just $Q \cdot (I_{Z_1})_3$, so this defines $Q \cup Z_1$
scheme-theoretically but is not saturated.  Notice that if we instead take
$d=6$, the ideal is both 6-regular and saturated, and the corresponding map is
injective.
\end{Ex}

\begin{Thm} \label{WLP delta 2}
Let $Z$ be a zero-dimensional subscheme of $\mathbb P^{n-1}$,
$n>3$, with \textup{WLP}. Suppose that
\[\Delta^2H(R/I_Z,d)=\Delta^2H(R/I_Z,d+1)=s\]
for $r_2(R/I_Z)>d>r_3(R/I_{Z})$, where $\Delta^2H(R/I_Z,\cdot)$ is
the second difference of Hilbert function of $R/I_Z$. Then
$\langle (I_Z)_{\leq d}\rangle$ is a saturated ideal defining a
two dimensional subscheme of degree $s$ in $\mathbb P^{n-1}$ and it
is $d$-regular.
\end{Thm}
\begin{proof}
With notation of Proposition \ref{Decreasing}, consider the exact
sequence (\ref{eq8}). By definition of WLP, we know that the map
\[(R/I_Z+(L_1))_d \stackrel{\times L_2}{\longrightarrow}
(R/I_Z+(L_1))_{d+1}\] is injective $d< r_2(R/I_Z)$ and
\[\Delta^2H(R/I_Z,t)=H(R/(I_Z+J),t)\]
for $t\leq r_2(R/I_Z)$ where $J=(L_1,L_2)$ for general linear
forms $L_1$ and $L_2$. Then the result follows from Theorem
\ref{r2 case}.
\end{proof}

\begin{Thm}\label{UPP_in_P3}
Let $Z$ be a set of points with UPP in $\mathbb P^{3}$. Let
$I_Z=(F_1,\ldots,F_m)$ be the defining ideal of $Z$, where
$\deg(F_i)=d_i$ and $d_1\leq d_2\leq \cdots \leq d_m$. Let $K =
(I_Z + (L_1,L_2))/(L_1,L_2)$ and let $S=R/(L_1,L_2)$. Then
\begin{equation}\label{eq7}
H(S/K,d)> H(S/K,d+1)
\end{equation}
for $d_2\leq d \leq r_2(R/I_Z)$ and $H(S/K,d) = 0$ for $d > r_2(R/I_Z)$.
\end{Thm}
\begin{proof}
First note that if $d \geq r_2(R/I_Z)$ then the assertions are just from the
definition of $r_2(R/I_Z)$.  If
$H(S/K,d)=H(S/K,d+1)$ for
$d<r_2(R/I_Z)$ then we know that
$(I_{Z})_d$ has a common factor $F$ by Corollary \ref{Common
factor}. Then, $F$ must be irreducible and $(I_Z)_{\leq d}=(F)$
since $Z$ has UPP (Lemma 4.4 in \cite{BGM}). This means that
$d_1\leq d<d_2$ and $\deg(F)=d_1$.
\end{proof}

\begin{Ex} \label{UPP vs WLP}
In Remark \ref{WLP connection} we discussed the connection between
the calculation
of the second reduction number for a zero-dimensional scheme and the
Weak Lefschetz
Property (WLP).  WLP  should be viewed as being a ``general'' property,
say for a fixed Hilbert function (as long as the Hilbert function
does not already
exclude the property).  But UPP is also a ``general'' property.  It
is easy to see
that a set of points can have WLP and not UPP (see \cite{HMNW} for
information on
the Hilbert function of a set of points with WLP).  In this example
we exhibit a
set of points, $\Gamma$, that has UPP but not WLP.

Let
$C$ be a smooth arithmetically Buchsbaum curve in
$\mathbb P^3$ whose deficiency module
$M(C)$  is one-dimensional in each of two consecutive degrees, and zero
elsewhere.  The Buchsbaum property means that all linear forms annihilate
$M(C)$.  The existence of a smooth curve in any even liaison class of curves in
$\mathbb P^3$ is proved in
\cite{rao}.  The precise shifts of $M(C)$ where smooth curves exist
were computed
in \cite{BM2} in a more general setting.  In this case, we may assume
that $M(C)$
is non-zero in degrees 3 and 4.

Let $L$ be a general linear form, defining a plane $H \subset
\mathbb P^3$.  Let $I_C$ be the saturated ideal of $C$.  We have
\[
J := I_C/(L \cdot I_C) \cong \frac{I_C + (L)}{(L)} \subset R/(L) := S
\]
and a short exact sequence of sheaves
\[
0 \rightarrow {\mathcal I}_C (-1) \stackrel{\times
L}{\longrightarrow} {\mathcal
I}_C \rightarrow {\mathcal I}_{C \cap H,H} \rightarrow 0.
\]
Taking cohomology on this last sequence and combining it with the isomorphisms
above, we obtain
\begin{equation} \label{isom1}
\frac{I_{C \cap H,H}}{J} \cong M(C)(-1)
\end{equation}
since $C$ is arithmetically Buchsbaum, where $I_{C\cap H,H}$ is the
saturation of
$\mathcal I_{C \cap H,H}$.

Now let $\Gamma$ consist of a general set of sufficiently many points on $C$.
The set $\Gamma$ has the Uniform Position Property since it consists
of a general
set of points on a smooth curve.  We want to show that $\Gamma$ does
not have the
Weak Lefschetz Property.

We have $(I_C)_t = (I_\Gamma)_t$ for all
$t \leq t_0$ for some $t_0$, which we may make as large as necessary
by choosing
sufficiently many points for $\Gamma$.
Since the points are general, $L$ does not vanish on any of them.
Then we have a
short exact sequence of sheaves
\[
0 \rightarrow {\mathcal I}_\Gamma (-1) \rightarrow {\mathcal I}_\Gamma
\rightarrow {\mathcal O}_H \rightarrow 0.
\]
Let $J' = I_\Gamma/(L \cdot I_\Gamma)$.  We obtain the isomorphism
\begin{equation} \label{isom2}
(S/J)_t = (S/J')_t \cong (\ker [ H^1_* ({\mathcal I}_\Gamma)(-1)
\stackrel{\times
L}{\longrightarrow} H^1_*({\mathcal I}_\Gamma)])_t
\end{equation}
for all $t \leq t_0$,
where $H^1_*$ refers to a direct sum over all twists. This isomorphism
commutes with multiplication by linear forms over this range of $t$.

Now consider the short exact sequence of sheaves
\[
0 \rightarrow {\mathcal I}_C \rightarrow {\mathcal I}_\Gamma \rightarrow
{\mathcal I}_{\Gamma | C} \rightarrow 0.
\]
Since $H^0(\mathcal I_{\Gamma|C}(t)) \cong H^0(\mathcal
O_C(tH-\Gamma))$ measures
the space  of hypersurface sections, up to linear equivalence,  of degree $t$
vanishing on
$\Gamma$ (but not on all of $C$), we may assume that
$h^0(\mathcal I_{\Gamma|C}(t)) = 0$ for $t \leq t_0$ (possibly changing $t_0$
slightly), again by taking sufficiently many points.  (In principle
$h^0(\mathcal
I_{\Gamma|C}(t))$ might be non-zero in degrees 3 and 4, coming from $M(C)$, but
then adding a hyperplane annihilates such a section, modulo $I_C$, which is
clearly nonsense.)  Hence we have
\[
0 \rightarrow H^1(\mathcal I_C(t)) \rightarrow H^1(\mathcal I_\Gamma(t))
\]
for $t \leq t_0$.  Therefore $M(C)$ is isomorphic to a submodule of
$H^1_*(\mathcal I_\Gamma)$.  Since
$M(C)$ is annihilated by all linear forms, we may invoke (\ref{isom2}) to
conclude that $M(C)$ is isomorphic to a submodule of $S/J'$.
In particular, since $M(C)$ occurs in degrees 3 and 4,
\begin{equation} \label{fact}
\hbox{$M(C)$ is isomorphic to a submodule of $S/J'$ occurring in
degrees 4 and 5,}
\end{equation}
thanks to
the shift in (\ref{isom2}).

   Let $L'$ be another general linear form, and consider the homomorphism
\[
(S/J')_4 \stackrel{\times L'}{\longrightarrow} (S/J')_5.
\]
We want to show that this is neither injective nor surjective.  The
fact that it
is not injective follows from (\ref{fact}).  For the non-surjectivity, consider
the commutative diagram
\[
\begin{array}{ccccccccccccccccccccccc}
&&&&&&0 \\
&&&&&& \downarrow \\ \\
0 & \rightarrow & \displaystyle \left ( \frac{I_{C\cap H,H}}{J} \right )_4 &
\rightarrow & \displaystyle \left ( \frac{S}{J} \right )_4 & \rightarrow &
\displaystyle \left ( \frac{S}{I_{C \cap H,H}} \right )_4 &
\rightarrow & 0 \\ \\
&& \phantom{\times 0} \downarrow { \times 0} && \phantom{\phi}
\downarrow \phi &&
\downarrow \\ \\
0 & \rightarrow & \displaystyle \left ( \frac{I_{C\cap H,H}}{J} \right )_5 &
\rightarrow & \displaystyle \left ( \frac{S}{J} \right )_5 & \rightarrow &
\displaystyle \left ( \frac{S}{I_{C \cap H,H}} \right )_5 &
\rightarrow & 0 \\ \\
&& \downarrow \\
&& A \\
&& \downarrow \\
&& 0
\end{array}
\]
The vertical arrows correspond to multiplication by $L'$.  Since $S/J =
S/J'$ in degrees 4 and 5 (and beyond), the middle vertical map is the
one that we
have to show is not surjective.  The fact that the leftmost vertical
map is zero
comes from (\ref{isom1}) and the Buchsbaum property.  Thanks to
(\ref{fact}), $A
\neq 0$.  The fact that the rightmost vertical map is injective comes from the
fact that $I_{C\cap H,H}$ is a saturated ideal.  Therefore we get the desired
non-surjectivity of the middle column from the Snake Lemma.
\end{Ex}

\begin{Ex} \label{hyper sect}
Uwe Nagel asked us whether our methods might be able to settle a 
question related to that of Example \ref{UPP vs WLP}.  That is, suppose that $C$
is a smooth curve over a field of characteristic zero.  Then it is well known
that the general hyperplane or hypersurface  section of $C$ has UPP.  Does
either of these  necessarily also have WLP?  When
$C \subset {\mathbb P}^3$, it is well known that the general hyperplane
section does have WLP, since all zero-dimensional schemes in ${\mathbb P}^2$ have
WLP
\cite{HMNW}.

A re-interpretation of Example \ref{UPP vs WLP} gives the surprising (to us)
answer ``no" to both questions.  Let $C$ be the smooth curve of degree 15 from
that example, and let $S$ be the cone over $C$ from a general point $P$ in
${\mathbb P}^4$.  After a change of variables, the generators of $I_C$, viewed
as polynomials in
$k[x_0,x_1,x_2,x_3,x_4]$, give the generators of $I_S$, and $S$ is smooth away
from $P$.  In particular, if $L=x_4$ is the linear form defining ${\mathbb P}^3$
in ${\mathbb P}^4$ (which holds without loss of generality after our change of
variables), then 
\begin{equation} \label{saturated}
I_C = [I_S + (L)]/(L)
\end{equation}  
where $I_C$ is the saturated ideal of $C$.  It follows from a
standard exact sequence that $H^1({\mathcal I}_S(t)) = 0$ for all $t$, since
(\ref{saturated}) implies that multiplication by a general linear form induces
an injection between any pair of consecutive components of
$H^1_*({\mathcal I}_S)$.
Consequently, for any homogeneous polynomial $F$ not vanishing on $S$ (which is
irreducible), 
$I_S + (F)$ is the saturated ideal of the hypersurface section of $S$ by $F$
(cf.\ \cite{migliore} Remark 2.1.3).

Now observe that in Example \ref{UPP vs WLP}, we could replace $\Gamma$ by a
general hypersurface section of sufficiently large degree.  Indeed, what is
important is that $(I_C)_t = (I_\Gamma)_t$ for all $t \leq t_0$, as describe
above.  This immediately shows that the general hypersurface section of a smooth
curve even in
${\mathbb P}^3$ does not necessarily have WLP.  But furthermore, let $L=x_4$ be
as above, and let $F$ be a general form  of sufficiently large degree in
$k[x_0,x_1,x_2,x_3,x_4]$, with $\bar F$ its restriction to
$k[x_0,x_1,x_2,x_3]$.  We have just observed that the zero-dimensional scheme
defined by $(I_C + \bar F)^{sat}$ fails to have WLP.  But 
\[
[I_C +(\bar F)]^{sat} = \left  [\frac{[I_S + (L)]}{(L)} + (\bar F) \right ]^{sat}
= \left [ \frac{[I_S +(F)] + (L)}{(L)} \right ]^{sat}.
\]
Since $F$ avoids $P$ for a general choice of $F$, $I_S +(F)$ is the saturated
ideal of a smooth curve, and the above equation shows that its hyperplane section
is the example we have already considered, which does not have WLP.
\end{Ex}

\begin{Remk}
The key idea in the preceding two examples is to use the structure of the
deficiency module to force the failure of WLP.  Consequently we do not know if
it is true that the general hyperplane or hypersurface section of a smooth
arithmetically Cohen-Macaulay curve necessarily has WLP.
\end{Remk}

Our last result gives further information about the growth of the Hilbert
function if we know UPP and WLP.

\begin{Coro}\label{UPP_WLP}
Let $Z$ be a set of points with UPP and WLP in $\mathbb P^{3}$.
Let $I_Z=(F_1,\ldots,F_m)$ be the defining ideal of $Z$, where
$\deg(F_i)=d_i$ and $d_1\leq d_2\leq \cdots \leq d_m$.
Then
\begin{equation}\label{eq7}
\Delta^2 H(R/I_Z,d)> \Delta^2 H(R/I_Z,d+1)
\end{equation}
for $d_2\leq d<r_2(R/I_Z)$.
\end{Coro}
\begin{proof}
The proof follows from Theorem \ref{UPP_in_P3} since
\[\Delta^2H(R/I_Z,t)=H(R/(I_Z+J),t)\] for $t\leq r_2(R/I_Z)$ where
$J=(L_1,L_2)$ for general linear forms $L_1$ and $L_2$ by WLP.
\end{proof}

\begin{Remk}
The result of Corollary \ref{UPP_WLP} holds only for points in ${\mathbb P}^3$. 
For ${\mathbb P}^n$, $n \geq 4$, we cannot say anything about the strictly
decreasing property of the second difference of the Hilbert function.  However,
it may well be that higher differences of the Hilbert function can be controlled
in higher projective spaces.
\end{Remk}

\bibliographystyle{amsalpha}

\end{document}